\documentclass[12pt]{article}
\usepackage{amssymb}
\usepackage{latexsym}
\usepackage{amsthm}  

\usepackage[dvips]{epsfig}

\newif\ifpdf
\ifx\pdfoutput\undefined
\pdffalse 
\else
\pdfoutput=1 
\pdftrue
\fi

\ifpdf
\usepackage[pdftex]{graphicx}
\else
\usepackage[dvips]{graphicx}
\fi

\newcommand{\N}{\mathbb N}
\newcommand{\Z}{ \mathbb Z}

\newcommand{\Sym}{\mathbb S}

\theoremstyle{plain}

\theoremstyle{remark}

\newenvironment{preuve}{\begin{trivlist}\item{\bf{Proof.}}}%
{\hfill\rule{2mm}{2mm}\end{trivlist}}
{\hfill\rule{2mm}{2mm}\end{trivlist}}

\newcommand{\ab}{\mbox{\Large\it a}}
\newcommand{\cat}{\mbox{\bf{c}}}

\newcommand{\abb}{{\mbox{\Large\it a}}^{-}}
\newcommand{\abbt}{{\mbox{\Large\it a}}^{\: \underline{\!\vartriangle\!}}}
\newcommand{\abt}{{\mbox{\Large\it a}}^{\vartriangle}}

\newcommand{\ap}[1][p]{\ab_{\mathrm{#1}}}
\newcommand{\apb}{\ap^{-}}
\newcommand{\apbt}{\ap^{\: \underline{\!\vartriangle\!}}}
\newcommand{\apt}{\ap^{\vartriangle}}

\newcommand{\tr}{\vartriangle}
\newcommand{\btr}{\: \underline{\!\vartriangle\!}}

\newcommand{\ars}{\mbox{$2$-trees}}
\newcommand{\ar}{\mbox{$2$-tree}}
\newcommand{\ie}{\textit{i.e.}\ }

\newtheorem{theo}{\bf Theorem}

\newtheorem{cor}{\bf Corollary}
\newtheorem{prop}{\bf Proposition}
\newtheorem{rem}{\bf Remark}

\begin{document}

\ifpdf
\DeclareGraphicsExtensions{.pdf, .jpg, .tif}
\else
\DeclareGraphicsExtensions{.eps, .jpg}
\fi

\title{A classification of plane and planar \ars{}}
\author{ Gilbert Labelle, C\'edric Lamathe, Pierre Leroux*\\
LaCIM, D\'epartement de Math\'ematiques, UQ\`AM}

\maketitle
\begin{abstract}
We present new functional equations for the species of plane and of planar (in the sense of Harary and Palmer, 1973) \ars{} and some associated pointed species. We then deduce the explicit molecular expansion of these species, \ie a classification of their structures according to their stabilizers. There result explicit formulas in terms of Catalan numbers for their associated generating series, including the asymmetry index series. This work is closely related to the enumeration of polyene hydrocarbons of molecular formula \(\mathrm{C}_n\mathrm{H}_{n+2}\).
\end{abstract}

\section{Introduction}
We define recursively the class $\ab$ of {\it 2-dimensional trees} (in brief {\it \ars{}}) as the smallest class of simple graphs such that 
\begin{itemize}
\item[1.] the single edge is in \ab,
\item[2.] if a simple graph $G$ has a vertex $x$ of degree 2 whose neighbors are adjacent and such that $G-x$ is in \ab, then $G$ is in \ab.
\end{itemize}
One can see that a \ar{} is essentially composed of triangles (complete graph on $3$ vertices) glued together along edges in a tree-like fashion.

\begin{figure}[ht]
 \centerline{\includegraphics[width=.6\textwidth]{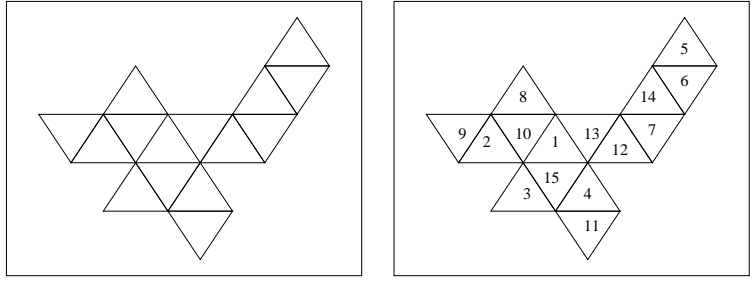}}
  \caption{An unlabelled plane \ar{} and one of its labellings}
  \label{ex1}
\end{figure}

\begin{figure}[ht]
 \centerline{\includegraphics[width=.6\textwidth]{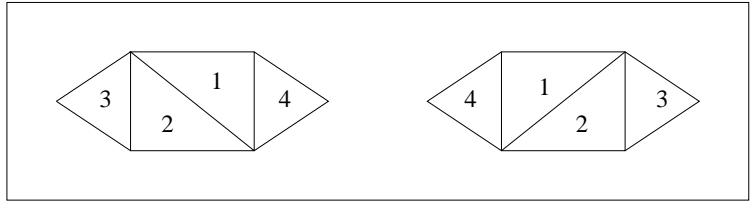}}
  \caption{Two different plane \ars{}, one planar \ar{}}
  \label{ex2}
\end{figure}

\begin{figure}[ht]
 \centerline{\includegraphics[width=.90\textwidth]{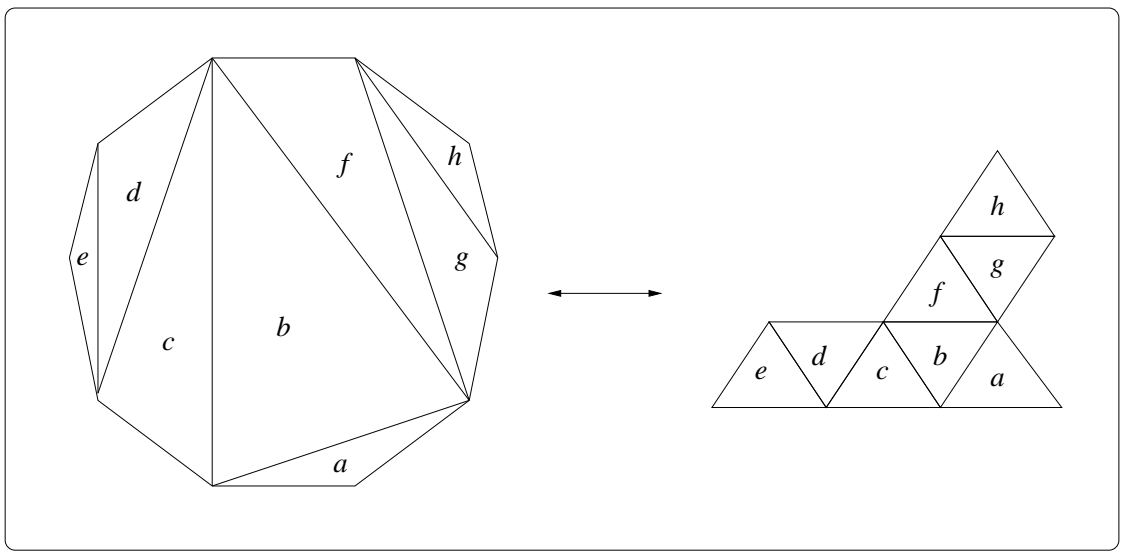}}
  \caption{Correspondence between triangulations of a polygon and plane \ars{}}
  \label{triang}
\end{figure}

Note that all \ars{} are planar simple graphs. However, by a {\it planar \ar}, we mean here a \ar{} admitting an embedding in the plane in such a way that all faces (except possibly the outer face) are triangles, and we call {\it plane \ar} a \ar{} equipped with such an embedding. This terminology agrees with Harary and Palmer~\cite{HP}. In Figure~\ref{triang}, we show a correspondence between plane \ars{} and (unrooted) triangulations of polygons in the plane which is also a correspondence between planar 2-trees and (unrooted) triangulations of polygons in space (no orientation), also known as triangulations of the disc, see \cite{BR} .
Figure~\ref{ex1} gives an example of an unlabelled and a triangle-labelled planar \ar{}, Figure~\ref{ex2} shows two different plane \ars{} which are in fact the same planar \ar{} since they are isomorphic simple graphs.
We point out the work of Palmer and Read, \cite{PR}, who enumerate plane embeddings of \ars{} without any condition on the faces, and which they also call plane \ars{}. 
 Planar \ars{} (in our sense) are closely related to acyclic polyene hydro-carbons of molecular formula \(\mathrm{C}_n\mathrm{H}_{n+2}\) (planar trees in the hexagonal lattice); see~\cite{CBBCL}.

We follow the approach of Fowler and al. in \cite{TF1,TF2} for general $\ars{}$. 
However, we go further here, giving explicitly the molecular expansion of plane and planar \ars{}, which could not be done in the general case. This is a stronger result than simple labelled and unlabelled enumeration since it gives a classification of the different structures according to stabilizers. For instance, it permits us to have an explicit enumeration of the symmetric and asymmetric parts of these species. Moreover, we obtain closed formulas for all coefficients appearing in these expansions.

To derive these results we use functional equations in the context of the combinatorial theory of species and deduce the molecular expansions and all the associated series. In the following, we label \ars{} at triangles and we denote by $X$ the species of singletons, \ie of simple triangles.
Recall that a combinatorial species is a class of finite labelled structures, closed under relabelling along bijections. To each species $F$ we associate series : $F(x)$, the exponential generating series of labelled structures; $\widetilde{F}(x)$, the ordinary generating series of unlabelled structures; ${\overline{F}}(x)$, the generating series of unlabelled asymmetric structures; $Z_F$ and $\Gamma_F$, the cycle and asymmetry index series. The usual shapes of these series for any species $F$ are as follows
\begin{equation}
F(x)=\sum_{n\geq 0}f_n{x^n\over n!},
\end{equation}
\begin{equation}
\widetilde{F}(x)=\sum_{n\geq 0}\widetilde{f}_nx^n,\quad \overline{F}(x)=\sum_{n\geq 0}\overline{f}_nx^n,
\end{equation}
\begin{equation}
Z_{F}(x_1,x_2,\ldots )=\sum_{n_1,n_2, \ldots}f_{n_1,n_2,\ldots }{x_1^{n_1}x_2^{n_2}\ldots \over 1^{n_1}n_1!2^{n_2}n_2!\ldots }\,,
\end{equation}
\begin{equation}
\Gamma_{F}(x_1,x_2,\ldots )=\sum_{n_1,n_2, \ldots}f_{n_1,n_2,\ldots }^*{x_1^{n_1}x_2^{n_2}\ldots \over 1^{n_1}n_1!2^{n_2}n_2!\ldots }\,,
\end{equation}
where $f_n$, $\widetilde{f}_n$ and $\overline{f}_n$ are the numbers of labelled, unlabelled and unlabelled asymmetric $F$-structures respectively, over an $n$-element set, and $f_{n_1,n_2,\ldots }$ is the number of $F$-structures left fixed under a given permutation of cycle type $1^{n_1}2^{n_2}\ldots$. For a definition of the asymmetry index series, see \cite{BLL}.\\

\begin{figure}[h]
 \centerline{\includegraphics[height=.25\textheight]{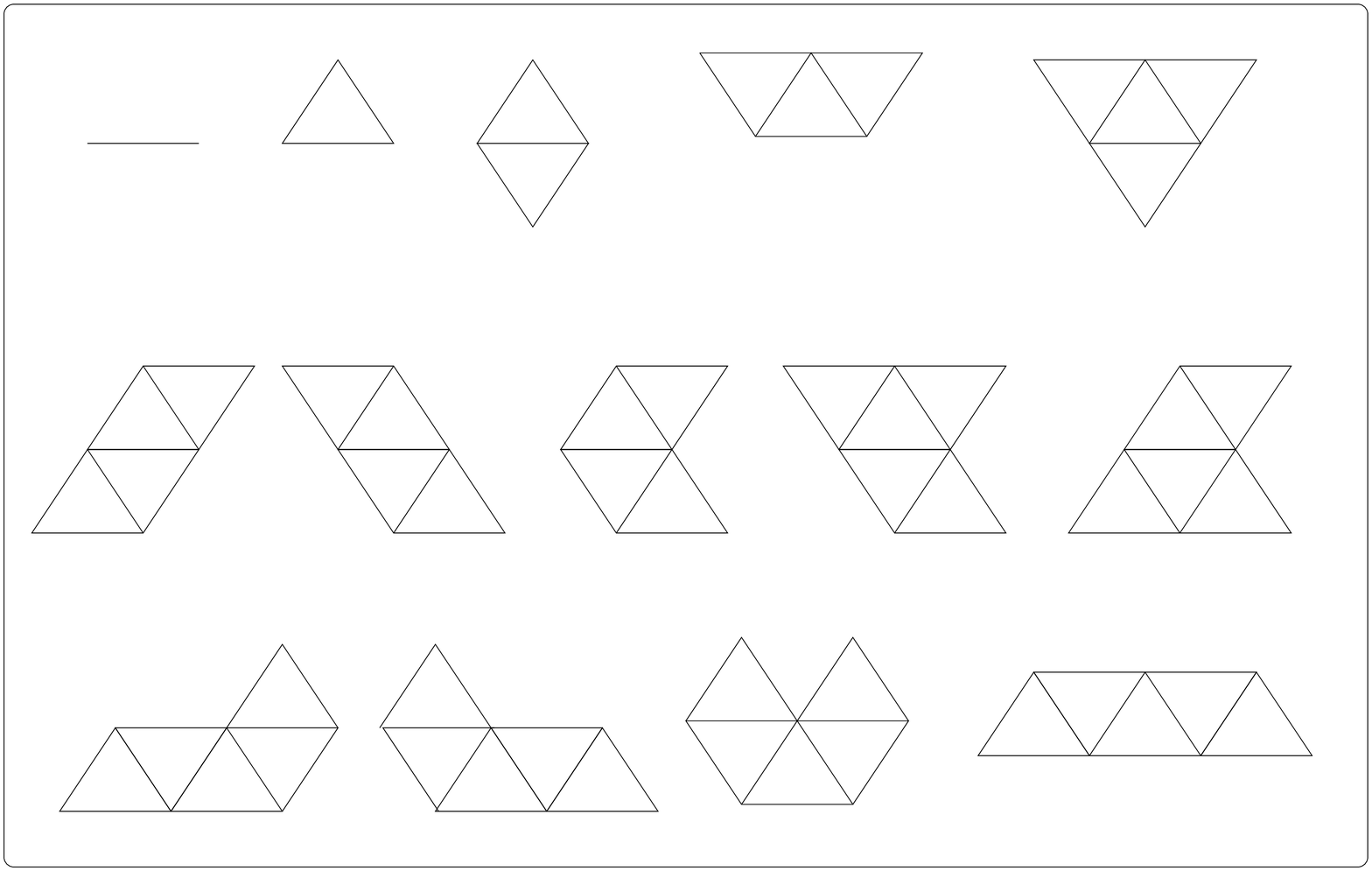}}
  \caption{First terms of the molecular expansion of the species $\ab_{\pi}$ of plane 2-trees}
  \label{decc1f}
\end{figure}

\begin{figure}[h]
 \centerline{\includegraphics[width=10cm]{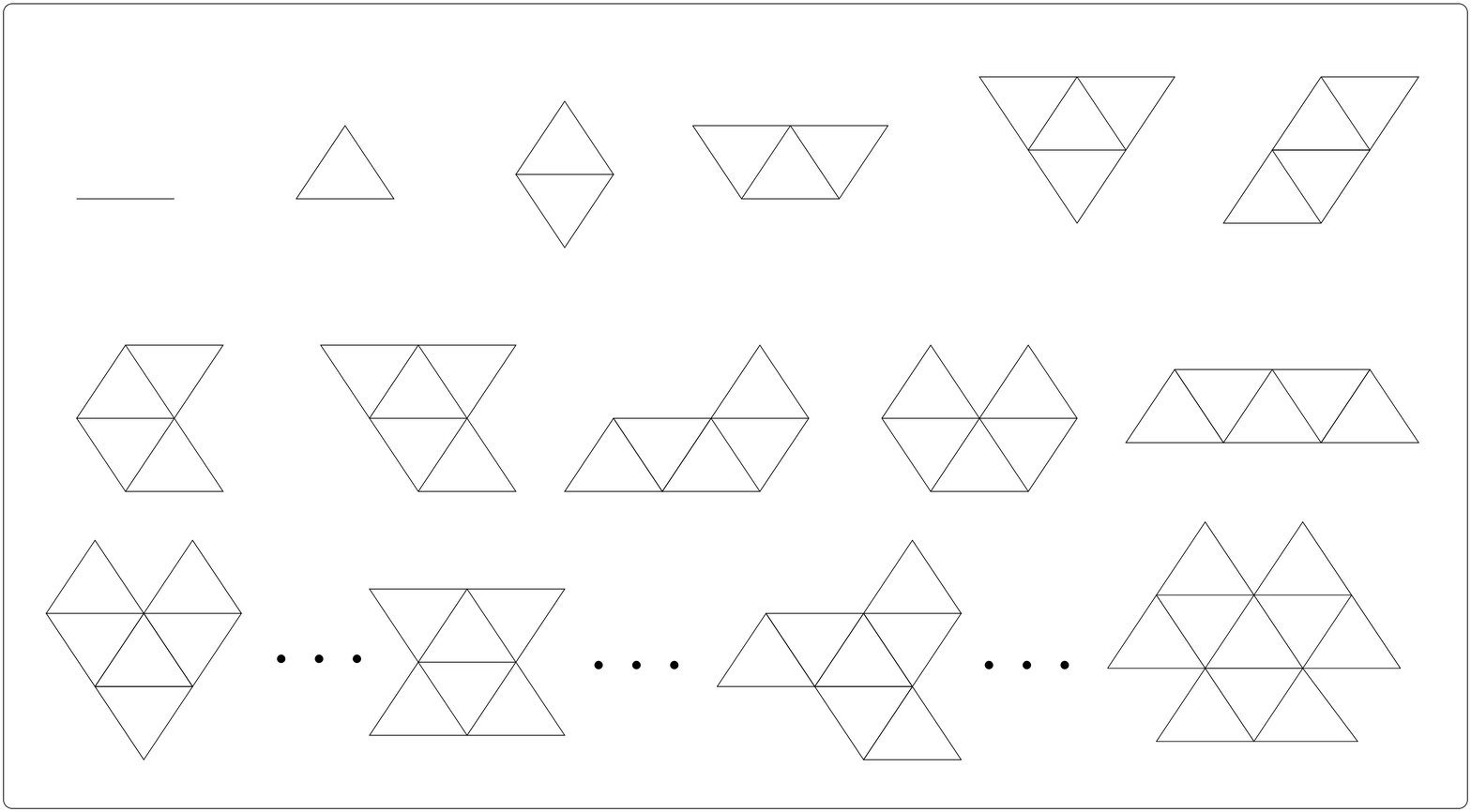}}
  \caption{First terms of the molecular expansion of the species $\ab_{\mathrm{p}}$ of planar 2-trees}
  \label{decc2f}
\end{figure}

To illustrate the notion of molecular expansion, we give here the first few terms of this decomposition for the species $\ab_{\pi}$ of plane \ars{} (Eq.\ (\ref{eq1}) and Figure~\ref{decc1f}) and $\ab_\mathrm{p}$ of planar \ars{} (Eq.\ (\ref{eq2}) and Figure~\ref{decc2f}). As usual, $E_n$ denotes the species of $n$-element sets and $C_3$, of 3-element (oriented) cycles. For complete explicit expansions see Theorem~\ref{Decmolplan} for plane 2-trees and Theorem~\ref{Decmolplanaire} for planar 2-trees.
\begin{equation}
\label{eq1}
\ab_{\pi}=\ab_{\pi}(X)=1+X+E_2(X)+X^3+XC_3(X)+2E_2(X^2)+X^4+6X^5+\cdots
\end{equation}
\pagebreak
$$
\ab_{\mathrm{p}}=\ab_{\mathrm{p}}(X)=1+X+E_2(X)+XE_2(X)+XE_3(X)+2E_2(X^2)+2X^5+2XE_2(X^2)
$$
\begin{equation}\label{eq2}
+X^2E_2(X^2)+\cdots +P_4^{bic}(X,X)+\cdots +XC_3(X^2)+\cdots +XP_6^{bic}(X,X)+\cdots\,.
\end{equation}

The expansion of $\ab_{\mathrm{p}}$ involves species $P_{4}^{bic}(X,Y)$ and $P_{6}^{bic}(X,Y)$ that are described in Section~$2$. They are two-sort variants of the species of $P_{2n}^{bic}$ introduced by J. Labelle in \cite{JL}.

\begin{figure}[ht]
 \centerline{\includegraphics[width=.55\textwidth]{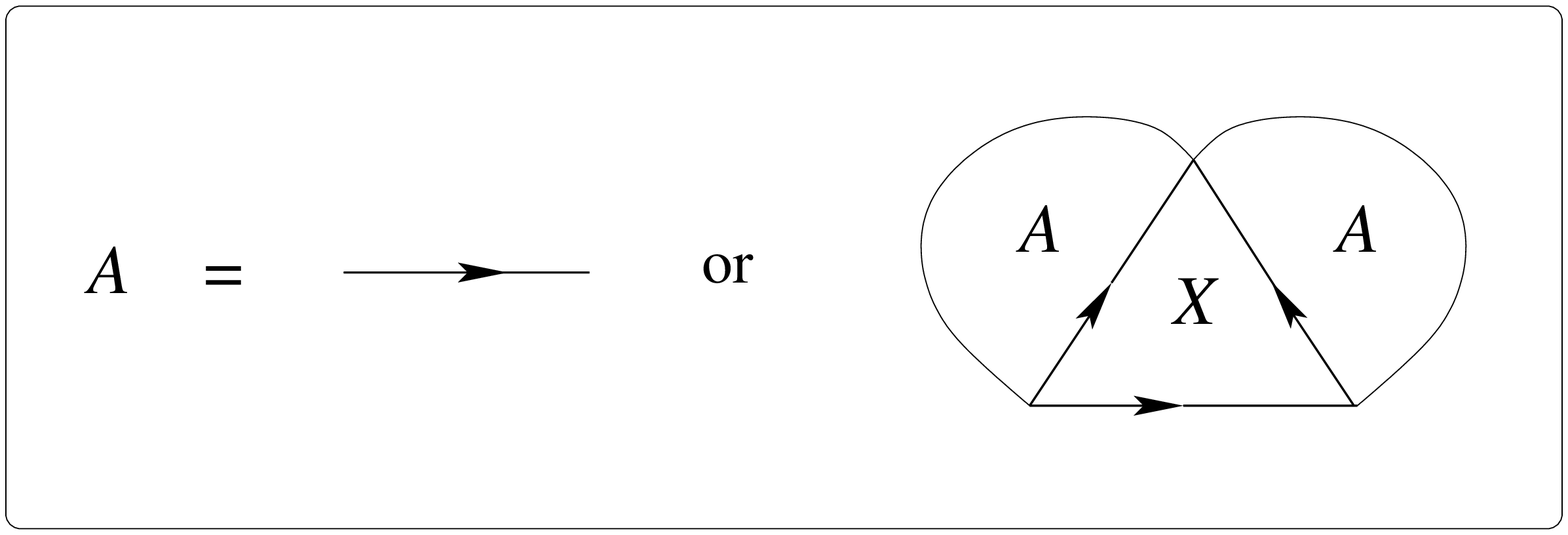}}
  \caption{$A=1+XA^2$}
  \label{ex3}
\end{figure}

In this paper, we call {\it degree} of an edge of a \ar{}, the number (less than or equal to 2) of triangles to which it belongs. Let us introduce the auxiliary species $A$ which can be defined as follows:
\begin{itemize}
\item[$\circ$] $A$ represents the species of plane \ars{} pointed at an external edge, \ie an edge of degree at most 1,
\item[$\circ$] $A$ is isomorphic to the species of planar \ars{} pointed at an external edge equipped with an orientation,
\item[$\circ$] $A$ is characterized by the functional equation 
\begin{equation}
A=1+XA^2\,,\label{A}
\end{equation}
illustrated in Figure~\ref{ex3}.
\end{itemize}
Note that the species $A$ can also be viewed as the species of rooted triangulations of polygons.
This species is fundamental for the following and we will use it several times. We can see that it is asymmetric, \ie the automorphism group of each of its structures is trivial; thus the molecular expansion and the associated series have the same coefficients in their expression. As expected, these coefficients are the Catalan numbers.
\begin{prop}\label{prop1}
The molecular expansion of the species $A=A(X)$ is
\begin{equation}
\label{eq11}
A(X)=\sum_{n\in \N}\cat_nX^n,
\end{equation}
where {$\cat_n={1\over n+1}{2n\choose n}$}  ({\it Catalan numbers}).
More generally, if {$A^k(X)=\sum_{n\in \N}\cat_{n}^{(k)}X^n$, $k\geq 1$}, then
\begin{eqnarray}
\cat_{n}^{(k)}&=&\sum_{i=0}^{\lfloor {k-1\over 2}\rfloor}(-1)^i{\scriptstyle{k-1-i}\choose \scriptstyle{i}}\cat_{n+k-1-i},\label{eq12}\\
&=&{k\over n}{2n-1+k\choose n-1}.
\label{eq13}
\end{eqnarray}%
\end{prop}
\begin{preuve}
The formula for $\cat_{n}$ follows directly from a simple application of the Lagrange inversion on the relation~(\ref{A}).
It can also be computed by expanding in series the algebraic solution
\(A(X)=(1-\sqrt{1-4X})/2X\) of~(\ref{A}).
%
For the \(\cat_{n}^{(k)}\), we work
with the unlabelled generating series. First, we remark that 
\begin{equation}
A^k(x)=\sum_{i=0}^{\lfloor {k-1\over 2}\rfloor}(-1)^i{\scriptstyle{k-1-i}\choose \scriptstyle{i}}{A(x)\over x^{k-1-i}}+\sum_{i=o}^{\lfloor {k-2\over 2}\rfloor }(-1)^{i+1}{\scriptstyle{k-2-i}\choose \scriptstyle{i}}{1\over x^{k-1-i}}, \label{merde}
\end{equation}
where $\lfloor \cdot \rfloor$ represents the floor function. This formula is easily shown by recurrence on $k$ distinguishing two cases depending on the parity of $k$ and using the fact that $A^2(x)={1\over x}(A(x)-1)$, which follows from (\ref{A}). Next, extracting the coefficient of $x^n$ in this expression gives the result. The second expression for $\cat_{n}^{(k)}$ is obtained by a simple application of the composite Lagrange inversion formula on equation (\ref{A}).
\end{preuve}

For instance, for $k$ from  $1$ up to $6$, we have 
\begin{eqnarray}
\cat_{n}^{(1)}&=&\cat_n\ =\ {1\over n}{2n\choose n-1}\,,\nonumber\\
\cat_{n}^{(2)}&=&\cat_{n+1}\ =\ {2\over n}{2n+1\choose n-1}\,,\nonumber\\
\cat_{n}^{(3)}&=&\cat_{n+2}-\cat_{n+1}\ =\ {3\over n}{2n+2\choose n-1}\,,\\
\cat_{n}^{(4)}&=&\cat_{n+3}-2\cat_{n+2}\ =\ {4\over n}{2n+3\choose n-1}\,,\nonumber\\
\cat_{n}^{(5)}&=&\cat_{n+4}-3\cat_{n+3}+\cat_{n+2}\ =\ {5\over n}{2n+4\choose n-1}\, ,\nonumber\\
\cat_{n}^{(6)}&=&\cat_{n+5}-4\cat_{n+4}+3\cat_{n+3}\ =\ {6\over n}{2n+5\choose n-1}\,.\nonumber
\end{eqnarray}

In order to lighten notations, we slighty extend the definition of the Catalan numbers as follows:
\begin{equation}
\cat_n={1\over n+1}{2n\choose n}\chi(n\in \N).\label{catalan}
\end{equation}
In other words, $\cat_n$ is the ususal Catalan number if $n$ is a nonnegative integer, and 0 otherwise.

We will use two dissymmetry formulas, analogous to the case of classical \ars{} (see Fowler and al. in \cite{TF1,TF2});
the same proof applies in the case of plane and planar \ars{} and is omitted.
\begin{theo}{{\textsc{Dissymmetry theorem for plane and planar \ars{}.}}}
The species $\ab_{\pi}$ of plane \ars{} and $\ab_{\mathrm{p}}$ of planar \ars{} satisfy the following isomorphisms of species
\begin{equation}
\ab_{\pi}^{-}+\ab_{\pi}^{\vartriangle}=\ab_{\pi}+\ab_{\pi}^{\: \underline{\!\vartriangle\!}},
\end{equation}
and
\begin{equation}
\ab_{\mathrm{p}}^{-}+\ab_{\mathrm{p}}^{\vartriangle}=\ab_{\mathrm{p}}+\ab_{\mathrm{p}}^{\: \underline{\!\vartriangle\!}},
\end{equation}
where the exponents $-$, $\vartriangle$ and $\underline{\!\vartriangle\!}$ represent the pointing of \ars{} at an edge (Figure~\ref{ex3a}a), at a triangle (Figure~\ref{ex3a}b) and at a triangle with one of its edges distinguished (Figure~\ref{ex3a}c).
\begin{figure}[h]
 \centerline{\includegraphics[width=.95\textwidth]{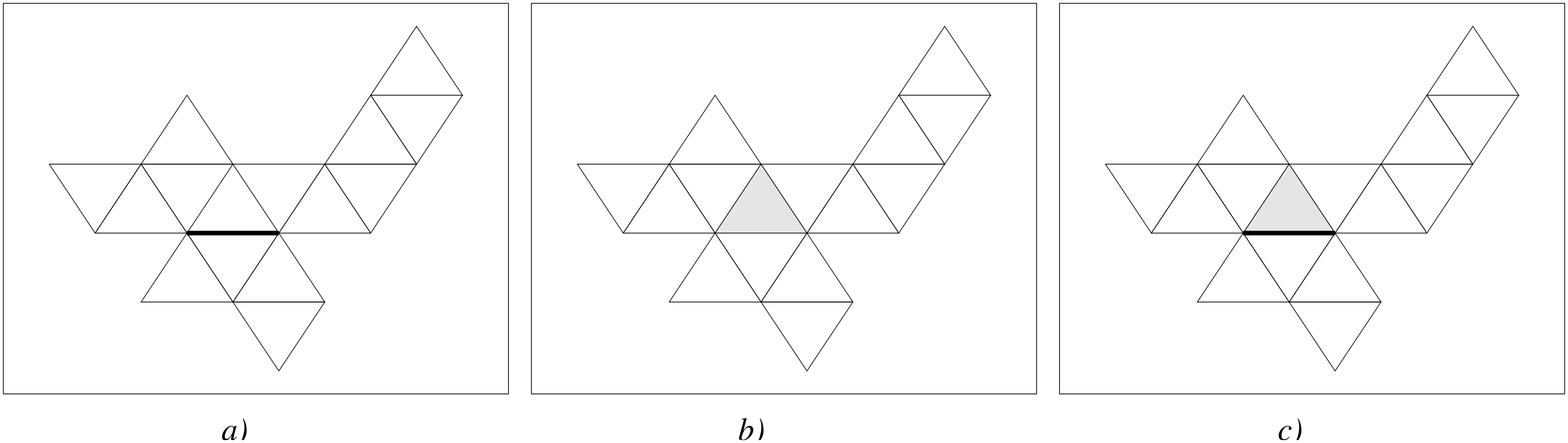}}
  \caption{Examples of the exponents: a) $-$, b) $\vartriangle$ and c) $\underline{\!\vartriangle\!}$}
  \label{ex3a}
\end{figure}
\end{theo}
The rest of the paper is organized as follows. In the next section, we introduce and study the auxiliary two-sort species $P_{4}^{\rm bic}(X,Y)$ and $P_6^{\rm bic}(X,Y)$ which are needed for the expression of the species $\apb$ and $\apt$ in terms of $A$. In Section 3, we give addition formulas for the substitution of an asymmetric species $Y=B(X)$ into the species $E_2(Y)$, $C_3(Y)$, $P_{4}^{\rm bic}(X,Y)$ and $P_6^{\rm bic}(X,Y)$. These results are put together in Section 4 to give the molecular expansion of the species $\ab_{\pi}$ and $\ap$. All the coefficients that occur in the expressions are given explicitly in terms of Catalan numbers. Finally, the labelled, unlabelled and asymmetric enumeration of plane and planar 2-trees is carried out in Section 5.

\section{The auxiliary molecular species $P_4^{\rm bic}(X,Y)$ and $P_6^{\rm bic}(X,Y)$}
This section is devoted to the study of some particular molecular species. A {\it molecular species} $M$ is a species having only one isomorphy type. In other words, any two $M$-structures are isomorphic. A molecular species is characterized by the fact that it is indecomposable under the combinatorial sum :
\begin{equation}
M\mbox{ is molecular}\quad\Leftrightarrow \quad (M=F+G\Rightarrow F=0\mbox{ or } G=0).
\end{equation}
It is often very useful to write a molecular species in the form
\begin{equation}
M={X^n\over H},
\end{equation}
where $X^n$ represents the species of lists of length $n$ and $H$ is a subgroup of the symmetric group $\Sym_n$. We write $H\leq \Sym_n$. In fact, $H$ is the stabilizer of some $M$-structure on $[n]=\{1,2\ldots n \}$ and $n$ is called the {\it degree} of the species $M$. Two molecular species of degree $n$, $X^n/H$ and $X^n/K$, are equal (\ie isomorphic as species) if and only if $H$ and $K$ are conjugate subgroups of $\Sym_n$.
 
Here are some examples of molecular species
\begin{itemize}
\item{} when $H=1$, then $X^n/1=X^n$,
\item{} when $H=\ <\rho>$, where $\rho$ is the circular permutation $\rho=(1,2,\ldots n)$, then $X^n/<\rho>\ =C_n$, the species of oriented cycles of length $n$,
\item{} if now the group $H$ is  $\Sym_n$, then we have $X^n/\Sym_n=E_n$, the species of sets of size $n$.
\end{itemize}
We denote by $\cal{M}$ the set of molecular species. We can see easily that the first elements of this set, up to degree 3, are
\begin{equation}
{\cal{M}}=\{1,X,X^2,E_2,X^3,XE_2,E_3,C_3(X),\ldots \}.
\end{equation}
Moreover, each species $F$ can be expressed as a (possibly infinite) linear combination with integer coefficients of molecular species as follows,
\begin{equation}
F=\sum_{M\in {\cal{M}}}f_MM,
\end{equation}
where $f_M\in \N$ represents the number of subspecies of $F$ isomorphic to $M$. This development is unique and it is called {\it molecular expansion} of the species $F$.

It is also possible to extend the notion of molecular species to the case of multi-sort species. For instance, for two-sort species, where $X$ and $Y$ represent  the two sorts, any molecular species can be written as 
\begin{equation}
M(X,Y)={X^nY^m\over H},\label{Mol}
\end{equation}
where $H\leq \Sym_n^X \times \Sym_m^Y$ is the stabilizer of an $M$-structure. Here, $\Sym_n^X$ represents the symmetric group of degree $n$ for the points of sort $X$.

\begin{figure}[h]
 \centerline{\includegraphics[width=.50\textwidth]{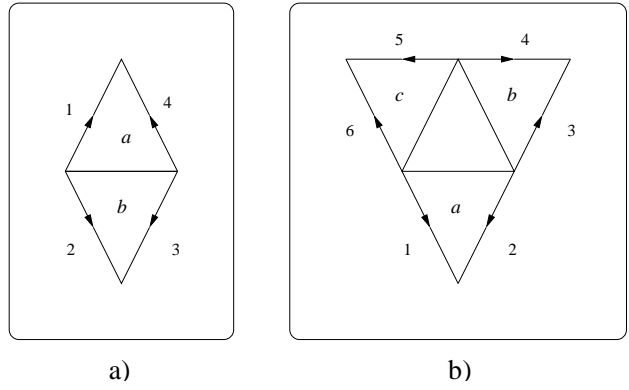}}
  \caption{Structures belonging to the species $Q(X,Y)$ and $S(X,Y)$}
  \label{f1}
\end{figure}
We can now introduce the auxiliary species $Q(X,Y)$ and  $S(X,Y)$ which will be important in our analysis of planar \ars{}. They can be defined by Figures~\ref{f1} a) and~\ref{f1} b) respectively,
where $X$ stands for the sort of triangles and $Y$, of directed edges.

\begin{figure}[h]
 \centerline{\includegraphics[width=16cm]{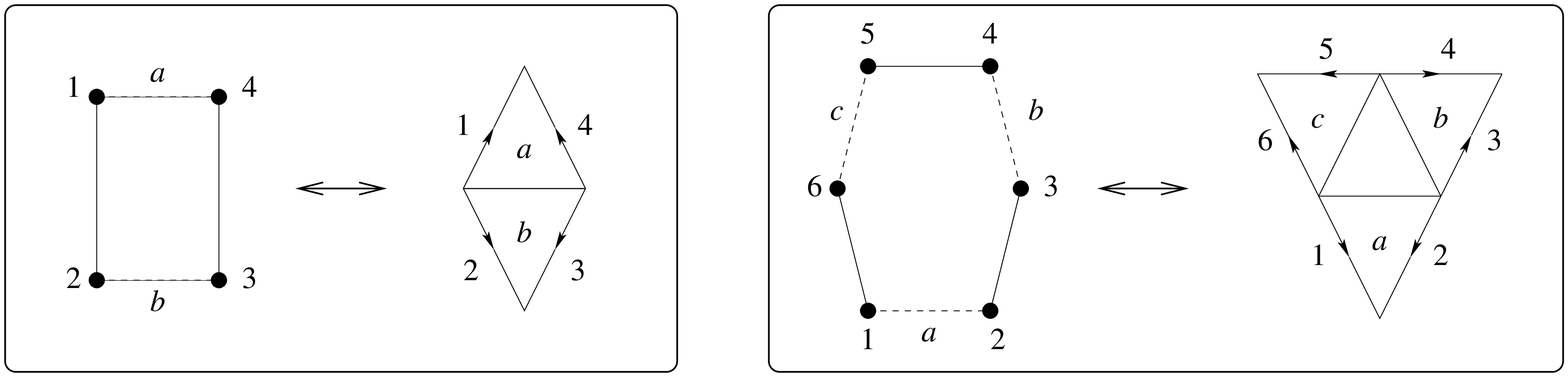}}
  \caption{$P_4^{\rm bic}(X,Y)$ and $P_6^{\rm bic}(X,Y)$}
  \label{f2}
\end{figure}

These two molecular species are related to known species:
 \begin{equation}
Q(X,Y)=P_4^{\rm bic}(X,Y),\quad  S(X,Y)=P_6^{\rm bic}(X,Y),\label{QS}
\end{equation}
where the species $P_n^{\rm bic}(X)$, for $n$ an even integer, represents the species of (vertex labelled) bicolored $n$-gons (see J. Labelle \cite{JL}). More precisely, the edges are colored with a set of two colors, $\{0,1\}$, in such a way that incident edges have different colors. We can then generalize to the two-sort species $P_n^{\rm bic}(X,Y)$ where $X$ represents the sort of edges of color $1$ (dotted lines) and $Y$ stands for the sort of vertices, as shown by Figure~\ref{f2} for $n=4$ and $n=6$. This Figure also establishes (\ref{QS}).

In order to completly describe the species $Q$ and $S$, we have to identify their stabilizers, and so we write them in the form (\ref{Mol}). We have 
\begin{equation}
P_4^{\rm bic}(X,Y)={X^2Y^4\over D_2},\quad P_6^{\rm bic}(X,Y)={X^3Y^6\over S_3}
\end{equation}
where the two groups $D_2$ and $S_3$ are characterized by their action on the labelled structures of Figure~\ref{f2} :
\begin{enumerate}
\item $D_2=<h,v> \leq \Sym_2^X\times\Sym_4^Y$, with  
$$
h=(a,b)(1,2)(3,4)\quad \mbox{and } \quad v=(a)(b)(1,4)(2,3).
$$
Note that $h^2=1,\ v^2=1,\ hv=vh,$ and $\ D_2 \stackrel{\sim}{=}\Z_2\times \Z_2$.
\item $S_3=<s,w> \leq \Sym_3^X\times\Sym_6^Y$, where 
$$
s=(a)(b,c)(1,2)(3,6)(4,5)\quad \mbox{and } \quad w=(a,b,c)(1,3,5)(2,4,6).
$$
Note that $s^2=1,\ w^3=1,\ sws=w^2,$ and $\ S_3 \stackrel{\sim}{=}\Sym_3$.
\end{enumerate}

Here are the formulas giving the cycle index series and the asymmetry index series of a molecular two-sort species.
\begin{theo}{\cite{BLL,K,LLP}}
\label{ZG}
Let $M(X,Y)=X^nY^m/H$ be a molecular species on two sorts, with $H\leq \Sym_n^X\times \Sym_m^Y$. Then, the cycle index series of $M$ is given by
\begin{equation}
Z_M(x_1,x_2, \ldots ;y_1,y_2,\ldots )={1\over |H|}\sum_{h\in H}x_1^{c_1(h)}x_2^{c_2(h)}\ldots   y_1^{d_1(h)}y_2^{d_2(h)}\ldots ,\label{ZZZ}
\end{equation} 
where $c_i(h)$ (resp. $d_i(h)$), for $i\geq 1$, denotes the number of cycles of lenght $i$ of the permutation on $X$-points (resp. $Y$-points) induced by the element $h\in H$. Furthermore, the asymmetry index series of $M$ is given by
\begin{equation}\label{gamma}
\Gamma_{M}( x_1,x_2, \ldots ;y_1,y_2,\ldots )={1\over |H|}\sum_{V\leq H}\mu(\{1\},V)x_1^{c_1(V)}x_2^{c_2(V)}\ldots  y_1^{d_1(V)}y_2^{d_2(V)}\ldots,    
\end{equation}
where the sum is taken over all subgroups $V$ of $H$, $\{1\}$ is the identity subgroup of $H$, $\mu(\{1\},V)$ denotes the value of the M\"obius function in the lattice of subgroup of $H$ and $c_i(V)$ (resp. $d_i(V)$), represents the number of orbits with $i$ elements of sort $X$ (resp. $Y$) with respect to the natural action of $V$ on $[n]$ (resp. $[m]$).
\end{theo}
\begin{prop}
The cycle index of the species $P_4^{\rm bic}(X,Y)$ and $P_6^{\rm bic}(X,Y)$ are given by
\begin{eqnarray}
Z_{P_4^{\rm bic}}(x_1,x_2,\ldots ;y_1,y_2,\ldots )&=&{1\over 4}(x_1^2y_1^4+2x_2y_2^2+x_1^2y_2^2),\\ \label{ZQ}
Z_{P_6^{\rm bic}}(x_1,x_2,\ldots ;y_1,y_2,\ldots )&=&{1\over 6}(x_1^3y_1^6+2x_3y_3^2+3x_1x_2y_2^3).\label{ZQa}
\end{eqnarray}
\end{prop}
\begin{preuve}
This is an easy exercise, using (\ref{ZZZ}) and writing explicitly the elements of the group $D_2$ and $S_3$ : \quad 
$\displaystyle{D_2=\{1,h,v,h\cdot v \}\quad \mbox{and}\quad S_3=\{1,s,\omega,\omega^2,s\cdot \omega,s\cdot \omega^2  \}}.$
\end{preuve}
\begin{prop}
The asymmetry index series of the two species $P_4^{\rm bic}(X,Y)$ and $P_6^{\rm bic}(X,Y)$ are given by
\begin{eqnarray}
\label{GQa}
\Gamma_{P_4^{\rm bic}}(x_1,x_2,\ldots ;y_1,y_2,\ldots )&=&{1\over 4}(x_1^2y_1^4-
x_1^2y_2^2-2x_2y_2^2+2x_2y_4),\\ \label{ZQ}
\Gamma_{P_6^{\rm bic}}(x_1,x_2,\ldots ;y_1,y_2,\ldots )&=&{1\over 6}(x_1^3y_1^6-x_3y_3^2-3x_1x_2y_2^3+3x_3y_6).\end{eqnarray}
\end{prop}
\begin{preuve}
It suffices to determine the lattice of subgroups of $D_2$ and $S_3$ and to apply (\ref{gamma}). Details are left to the reader.
\end{preuve}

The cycle index series of a species encompasses the two other classical enumerative series, namely the exponential generating function of labelled structures and the ordinary generating function of unlabelled structures. In a similar way, the asymmetry index series contains other series as specializations, in particular the asymmetry generating series. For the two-sort case, these series are related as follows :
\begin{theo}{(\cite{BLL}).}
For any two-sort species $F$, we have
\begin{eqnarray}
F(x,y)&=&Z_{F}(x,0,\ldots ;y,0,\ldots)\ =\ \Gamma_F(x,0,\ldots ;y,0,\ldots),\\
\widetilde{F}(x,y)&=&Z_F(x,x^2,\ldots ;y,y^2, \ldots),\\
\overline{F}(x,y)&=&\Gamma_F(x,x^2,\ldots ;y,y^2,\ldots).
\end{eqnarray}
\end{theo}
We then confirm the expressions of the generating series of the species $P_4^{\rm bic}(X,Y)$  and $P_6^{\rm bic}(X,Y)$.
\begin{rem}
We have 
\begin{eqnarray}
P_4^{\rm bic}(x,y)&=&{1\over 4}x^2y^4,\quad \widetilde{P}_4^{\rm bic}(x,y)=x^2y^4,\quad \overline{P}_4^{\rm bic}(x,y)=0,\\
P_6^{\rm bic}(x,y)&=&{1\over 6}x^3y^6,\quad \widetilde{P}_6^{\rm bic}(x,y)=x^3y^6,\quad \overline{P}_6^{\rm bic}(x,y)=0. 
\end{eqnarray}
\end{rem}
The fact that $\overline{P}_4^{\rm bic}(x,y)$ and $\overline{P}_6^{\rm bic}(x,y)$ equals $0$, means that these two species are purely symmetric, \ie, their asymmetric part is reduced to the empty set.\\

Note that if we put $Y:=X^k$, for $k\geq 1$, in the species $P_4^{\rm bic}(X,Y)$ and $P_6^{\rm bic}(X,Y)$, the resulting  one-sort species are molecular. Indeed, the substitution of a molecular species in another one remains molecular. These two species $P_4^{\rm bic}(X,X^k)$ and $P_6^{\rm bic}(X,X^k)$, for $k\geq 1$, will be essential in order to obtain the molecular expansion of planar 2-trees. Besides, we remark the fact that 
\begin{equation}
P_4^{\rm bic}(X,1)=E_2(X),\quad  P_6^{\rm bic}(X,1)=E_3(X),
\end{equation}
since, in Figure 8, setting $Y=1$ corresponds to unlabelling the directed edges.
 
To end this section, let us give the derivative of the two-sort species $P_4^{\rm bic}(X,Y)$ and $P_6^{\rm bic}(X,Y)$.
\begin{prop} \label{derivee}
The partial derivatives of $P_4^{\rm bic}(X,Y)$ and $P_6^{\rm bic}(X,Y)$ are given by
\begin{eqnarray}
{\partial \over \partial X}P_4^{\rm bic}(X,Y)=XE_2(Y^2),\quad   {\partial \over \partial Y}P_4^{\rm bic}(X,Y)=X^2Y^3,\label{der1}\\
{\partial \over \partial X}P_6^{\rm bic}(X,Y)= E_2(XY^3),\quad   {\partial \over \partial Y}P_6^{\rm bic}(X,Y)=X^3Y^5.\label{der2}
\end{eqnarray}
\end{prop}
\begin{preuve}
Let $F(X,Y)$ be a two-sort species and $U$ and $V$ be two sets representing the two sorts. Then, the partial derivatives, with respect to $X$ and $Y$ are defined by
$$
{\partial F\over \partial X}[U,V] = F[U+\{*\},V],\quad 
{\partial F\over \partial Y}[U,V] = F[U,V+\{*\}],
$$
where $*$ is a supplementary element which is used in the construction of the $F$-structures. From this definition, it is easy to obtain (\ref{der1}) et (\ref{der2}).
\end{preuve}

\section{Addition formulas}
In this section, we prove some addition formulas which will be necessary to obtain the explicit molecular expansions for plane and planar 2-trees.
\begin{prop}
\label{lem3}
Let $B$ be an asymmetric species whose molecular expansion is given by 
$$
{B(X)=\sum_{k\geq 0}b_kX^k}\,.
$$
Then, we have the following addition formulas relative to the species $E_2$ of two-element sets and $C_3$ of oriented 3-cycles :
\begin{eqnarray}
E_2(B(X))&=&\sum_{k\geq 1}b_kE_2(X^k)+\sum_{k\geq 0}\alpha_kX^k, \label{E2}\\
C_3(B(X))&=&\sum_{k\geq 1}b_kC_3(X^k)+\sum_{k\geq 0}\beta_kX^k, \label{C3}
\end{eqnarray}
with
\begin{eqnarray}
\lefteqn{%
\alpha_0= \frac{1}{2}(b_0^2+b_0),
\quad
 \beta_0 = \frac{1}{3}(b_0^3+2b_0)\,,} &&\label{alpha0}\\
\alpha_k&=&{1\over 2}\sum_{i+j=k}b_ib_j-{1\over 2}\chi(2|k)b_{{k\over 2}},
\quad k\geq1\,,\label{alpha}\\
\beta_k&=&{1\over 3}\sum_{l+m+n=k}b_lb_mb_n-{1\over 3}\chi(3|k)b_{{k\over 3}},
\quad k\geq1\,,\label{beta}
\end{eqnarray}
where, for $a,b \in \N$, $\chi(a|b)=1$, if $a$ divides $b$, and $0$, otherwise.
\end{prop}

\begin{preuve}
First note that for any species $F$, the constant (\ie of degree 0) term $F(b_0)$ of $F(B)$ is given by $Z_F(b_0,b_0,\ldots )$, in virtue of Polya's theorem. This yields (\ref{alpha0}). An analysis of the different shapes of molecular species which can arise in $E_2(B)$, permits us to write the following relation 
\begin{equation}
E_2(B)=\sum_{k\geq 1}\gamma_kE_2(X^k)+\sum_{k\geq 0}\alpha_kX^k. \label{q}
\end{equation}
We now have to compute $\alpha_k$ and $\gamma_k$, for all $k\geq 1$. Note that we can order, in the species $B$, the $b_k$ copies of the molecule $X^k$, for each $k\geq 1$. Then, to obtain an $E_2(X^k)$-structure from $E_2(B)$, we must take twice the same copy of $X^k$ among the $b_k$ available; otherwise the pair of $B$-structures will be asymmetric. Hence $\gamma_k=b_k$, for all $k\geq 1$. In order to compute $\alpha_k$, we could perform a direct enumeration. However, we introduce a different method which will prove very useful in other situations. Differentiating the two members of (\ref{q}), we get
$$
BB^{\prime}=\sum_{k\geq 1}kb_kX^{2k-1}+\sum_{k\geq 1}k\alpha_kX^{k-1}.
$$
Integrating back this last relation, in the realm of formal power series in $X$, leads us to
$$
{1\over 2}B^2={1\over 2}\sum_{k\geq 1}b_kX^{2k}+\sum_{k\geq 0}\alpha_kX^k+\mbox{const}\,.
$$
Identifying coefficients of $X^n$ in both sides of the last equality gives us the relation (\ref{alpha}).
To obtain (\ref{beta}), we first write
\begin{equation}
C_3(B)=\sum_{k\geq 1}\delta_kC_3(X^k)+\sum_{k\geq 0}\beta_kX^k. \label{qq}
\end{equation}
The same argument as used above implies $\delta_k=b_k$, $k\geq 1$, and the same technique of differentiating-integrating equation (\ref{qq}) gives the announced formula for $\beta_k$. In the process, we use the fact that
$$
(C_3(B))^{\prime}=L_2(B)B^{\prime}=B^2B^{\prime}
$$
where $L_2$ represents the species of two-element lists.
\end{preuve}

As a  particular case, we have
\begin{eqnarray}
E_2(1+X)&=&1+X+E_2(X),\\
C_3(1+X)&=&1+X+X^2+C_3(X).
\end{eqnarray}

When $B=A$, formulas (\ref{alpha0})--(\ref{beta}) take a simpler form because of the convolutive properties of Catalan numbers, as seen in Proposition~\ref{prop1}. For this case, the coefficients $\alpha_k$ and $\beta_k$ are given by \(\alpha_0=\beta_0=1\) and, for \(k\geq1\),
\begin{eqnarray}
\alpha_k&=&{1\over 2}(\cat_{k+1}-\cat_{{k\over 2}}),\\
\beta_k&=&{1\over 3}(\cat_{k+2}-\cat_{k+1}-\cat_{{k\over 3}}).
\end{eqnarray}
We now give the main result of this section, addition formulas for the species $P_4^{\rm bic}(X,Y)$ and  $P_6^{\rm bic}(X,Y)$. Let $b_k^{(n)}$ denotes the coefficient of $X^k$ in the species $B^n(X)$, with the convention that $b_x^{(n)}=0$ if the index $x$ is fractional, for all $n,k\geq 1$. 
\begin{theo}
Let $B$ be an asymmetric species whose molecular expansion is given by $$\displaystyle{B(X)=\sum_{n\geq 0}b_kX^k}.$$ Then, 
\begin{equation}\label{addP4}
P_4^{\rm bic}(X,B)={{\displaystyle{\sum_{k\geq 3}}}a_k^{'}X^k+{\displaystyle{\sum_{k\geq 2}}}a_k^{''}E_2(X^k)+{\displaystyle{\sum_{k\geq 1}}}a_k^{'''}X^2E_2(X^k)+{\displaystyle{\sum_{k\geq 0}}}a_k^{iv}P_4^{\rm bic}(X,X^k)},
\end{equation}
where
\begin{eqnarray}
a_k^{'}&=&{1\over 4}b_{k-2}^{(4)}-{3\over 4}b_{k-2\over 2}^{(2)}+{1\over 2}b_{k-2\over 4},\label{P41}\\
a_k^{''}&=&b_{k-1}^{(2)}-b_{k-1\over 2},\label{P42}\\
a_k^{'''}&=&{1\over 2}(b_{k}^{(2)}-b_{k\over 2}),\label{P43}\\
a_k^{iv}&=&b_k.\label{P44}
\end{eqnarray}
\end{theo}
\begin{preuve}
We proceed in a similar way as in Proposition~\ref{lem3}, beginning with an analysis of the different symmetries which can appear in structures belonging to the species $P_4^{\rm bic}(X,B(X))$. This permits us to write (\ref{addP4})
where all coefficients have to be determined. We first note that $a_k^{iv}=\cat_k$ since the only way to build a $P_4^{\rm bic}(X,X^k)$-structure from the species $P_4^{\rm bic}(X,B)$ is to take four times the same copy of the molecule $X^k$ among the $b_k$ available copies. This gives (\ref{P44}). Next, we consider $E_2(X^k)$-structures. In order to obtain such a structure from the species $P_4^{\rm bic}(X,B)$, we can take two non isomorphic $X^{{k-1\over 2}}$-structures $\alpha$ and $\beta$ from the species $B$, and put them in the two differents ways shown in Figure \ref{12} a) and \ref{12} b). This contributes for a term of
$$
2\sum_{l}{b_l\choose 2}E_2(X^{2l+1}),
$$
remembering that the two internal triangles also contribute for one $X$ each. 
We can also take an $X^i$-structure $\alpha$ and an $X^j$-structure $\beta$ such that $i+j=k-1$ and $i\neq j$, and put them in the two different configurations drawn in Figure~\ref{12} $a)$ and $b)$.
\begin{figure}[ht]
 \centerline{\includegraphics[width=.60\textwidth]{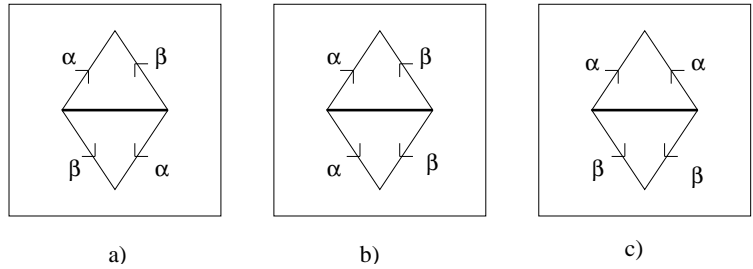}}
 \caption{Symmetries of order 2 in $P_4^{\rm bic}(X,B)$}
  \label{12}
\end{figure}
In the molecular expansion of the species $Q(X,B)$ this stands for
$$
2\sum_{i+j=k-1\atop i<j}b_ib_jE_2(X^k).
$$
It leads to (\ref{P42}), \ie
$$
a_k^{''}=2\sum_{i+j=k-1\atop i<j}b_ib_j+{b_{k-1\over 2}\choose 2}=b_{k-1}^{(2)}-b_{k-1\over 2}.
$$
Let us now turn to the coefficient $a_k^{'''}$ of $X^2E_2(X^k)$ in the relation (\ref{addP4}). The configurations belonging to an $X^2E_2(X^k)$ are shown in Figure~\ref{12} $c)$. We then have 
$$
a_k^{'''}=\sum_{i+j=k\atop i<j}b_ib_j+{b_{k\over 2}\choose 2}={1\over 2}(b_k^{(2)}-b_{k\over 2})
$$
types of $X^2E_2(X^k)$-structures. It remains to determine the asymmetric part of the species $Q(X,B)$, \ie the coefficient $a_k^{'}$ of $X^k$ in the molecular expansion (\ref{addP4}), for all $k$. To find it, we differentiate the relation (\ref{addP4}) and  we identify the coefficient of $X^k$ in each side. It gives the expression (\ref{P41}), which completes the proof. Note that we use the combinatorial derivative of a composite species $F(X,B(X))$. As in calculus, we have 
\begin{equation}
(F(X,B(X)))^{\prime}={\partial F(X,Y)\over \partial X}|_{Y:=B}+{\partial F(X,Y)\over \partial Y}|_{Y:=B}\cdot B^{\prime},
\end{equation}
and we can use Proposition~\ref{derivee}.
\end{preuve}
\begin{rem}
We can perform a precise classification separating rotational and reflectional symmetries. Indeed, the symmetries illustrated by Figure \ref{12} are rotational for the case a), vertically reflectional for case b) and horizontally reflectional for c).
\end{rem}
Remark also that we could obtain the expression of $a_k^{'''}$ by identifying the coefficient of $XE_2(X^k)$ after deriving (\ref{addP4}).
\begin{theo}
For all asymmetric species $B$ whose molecular expansion is 
$$
{B(X)=\sum_{k\geq 0}b_kX^k} \,,
$$ 
we have 
\begin{equation}
\label{nowhere}
P_6^{\rm bic}(X,B)={{\displaystyle{\sum_{k\geq 4}}}d_k^{'}X^k+{\displaystyle{\sum_{k\geq 2}}}d_k^{''}XE_2(X^k)+{\displaystyle{\sum_{k\geq 2}}}d_k^{'''}C_3(X^k)+{\displaystyle{\sum_{k\geq 0}}}d_k^{iv}P_6^{\rm bic}(X,X^k)},
\end{equation}
where
\begin{eqnarray}
d_k^{'}&=&{1\over 6}b_{k-3}^{(6)}-{1\over 2}b_{k-3\over 2}^{(3)}+{1\over 3}b_{k-3\over 3}^{(2)}+{2\over 3}b_{k-3\over 6},\\
d_k^{''}&=&b_{k-1}^{(3)}-b_{{k-1\over 3}},\\
d_k^{'''}&=&{1\over 2}(b_{k-1}^{(2)}-b_{k-1\over 2}),\\
d_k^{iv}&=&b_k,
\end{eqnarray}
where $b_k^{(n)}$ represents the coefficient of $X^k$ in $B^n(X)$.
\end{theo}
\begin{preuve}
A precise analysis of the different symmetries arising in the species $P_6^{\rm bic}(X,B)$ permit us to write the expansion (\ref{nowhere}). We then compute all coefficients of this expression by the same method as for the species $P_4^{\rm bic}(X,B)$.
\end{preuve}

When we put $B=A$ in the two previous theorems, the coefficients appearing in the molecular expansions of the species $P_4^{\rm bic}$ and $P_6^{\rm bic}$ are simpler. In fact, by Proposition~\ref{prop1} we get the following expressions for $a_k^{i}$ and $d_k^i$, for $i\in\{\prime,\prime \prime,\prime \prime \prime,iv\}$ 

\begin{eqnarray}
a_k^{'}&=&{{1\over 4}}\cat_{k+1}-{{1\over 2}}\cat_{k}-{{3\over 4}}\cat_{k\over 2}+{{1\over 2}}\cat_{k-2\over 4},\nonumber\\
a_k^{''}&=&\cat_k-\cat_{k-1\over 2},\\
a_k^{'''}&=&\frac{1}{2}(\cat_{k+1}-\cat_{k\over 2}),\nonumber\\
a_k^{iv}&=&\cat_k, \nonumber\\
&&\nonumber\\
d_k^{'}&=&{{1\over 6}}\cat_{k+2}-{2\over 3}\cat_{k+1}+{1\over 2}\cat_k-{1\over 2}\cat_{k+1\over 2}+{1\over 2}\cat_{k-1\over 2}-{1\over 6}\cat_{k\over 3}+{1\over2}\cat_{k-3\over 6},\nonumber\\
d_k^{''}&=&\cat_{k+1}-\cat_k-\cat_{\frac{k-1}{3}},\\
d_k^{'''}&=&\frac{1}{2}(\cat_{k}-\cat_{k-1\over 2}),\nonumber\\
d_k^{iv}&=&\cat_k\nonumber.
\end{eqnarray}                         
%
%
%
%
%
%
%
%
%
%
%
%
%
\section{Molecular expansion of plane and planar \ars{}}
In this part, we use the dissymmetry theorem and the results of the previous section to obtain an explicit form for the molecular expansion of the species of plane 2-trees and of planar 2-trees.
\subsection{Plane $2$-trees}
\begin{figure}[ht]
 \centerline{\includegraphics[height=.175\textheight]{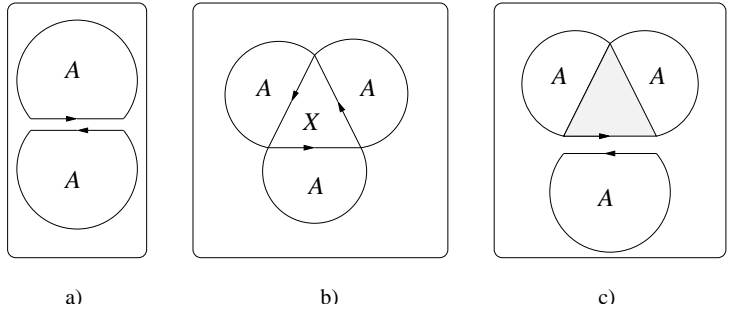}}
  \caption{The species $E_2(A)$, $XC_3(A)$ and $A_+\cdot A$}
  \label{th1}
\end{figure}

Recall that plane 2-trees are 2-trees that are embedded (drawn) in the plane in such a way that all internal faces are triangles. The dissymetry theorem gives an expression for the species $\ab_{\pi}$ in terms of the pointed species $\abb_{\pi}$, $\abt_{\pi}$ and $\abbt_{\pi}$, namely 
\begin{equation}
\label{DTH}
\ab=\abb_{\pi}+\abt_{\pi}-\abbt_{\pi}.
\end{equation}
Here, we can use the orientation of the plane to obtain simple expressions for the pointed species as function of the species $A$ defined in the introduction, as shown in Figure~\ref{th1} :
\begin{theo}
The species arising in the dissymmmetry theorem for plane 2-trees satisfy
\begin{eqnarray}
\ab_{\pi}^{-}&=&E_{2}(A), \label{t}\\
\ab_{\pi }^{\vartriangle}&=&XC_{3}(A),\label{tt}\\
\ab_{\pi}^{\:\underline{\!\vartriangle\!}}&=&A_+\cdot A,\label{ttt}
\end{eqnarray}
where $A_+=A-1$.
\end{theo}
%

Using the expansion formulas for $E_2(A)$ and $C_3(A)$, given in Section 3, we can now compute       the molecular expansion of the species $\ab_{\pi}$.
\begin{theo}\label{Decmolplan}
The molecular expansion of the species $\ab_{\pi}$ of plane \ars{} is given by
\begin{equation}
\ab_{\pi}=\ab_{\pi}(X)=1+X+\sum_{k\geq 2}b_kX^k+\sum_{k\geq 1}c_kE_2(X^k)+\sum_{k\geq 1}d_kXC_3(X^k), \label{DEC1}
\end{equation}
where 
\begin{eqnarray}
b_k&=&{2\over 3}\cat_k-{1\over 6}\cat_{k+1}-{1\over 2}\cat_{{k\over 2}}-{1\over 3}\cat_{{k-1\over 3}},\label{coef1}\\
c_k&=&d_k\quad =\quad \cat_k, \label{coef2}
\end{eqnarray}
where $X^k$ represents the species of $k$-lists of triangles and $\cat_k$ are the usual Catalan numbers with the convention that $\cat_r=0$ if $r$ is not an integer; see (\ref{catalan}). 
\end{theo}


To conclude this section we write the asymmetric part, in the sense of G.\ Labelle \cite{GL}, of the species of plane \ars{} : 
\begin{equation}
\overline{\ab}_{\pi}(X)=1+X+\sum_{k\geq 2}b_kX^k,
\end{equation}
where $b_k$, for $k\in \N$, is given by the formula (\ref{coef1}). The species $\overline{\ab}_{\pi}$ is not to be confused with the pointed species $\ab_{\pi}^{-}$.

\subsection{Planar \ars{}}
\begin{figure}[ht]
 \centerline{\includegraphics[height=.125\textheight]{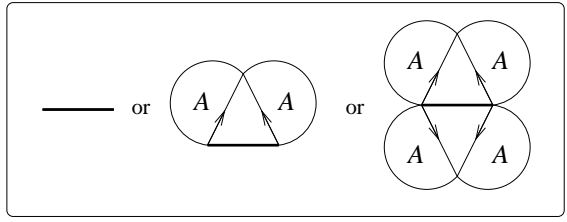}}
  \caption{The species $\apb$}
  \label{new1}
\end{figure}

\begin{figure}[ht]
 \centerline{\includegraphics[height=.18\textheight]{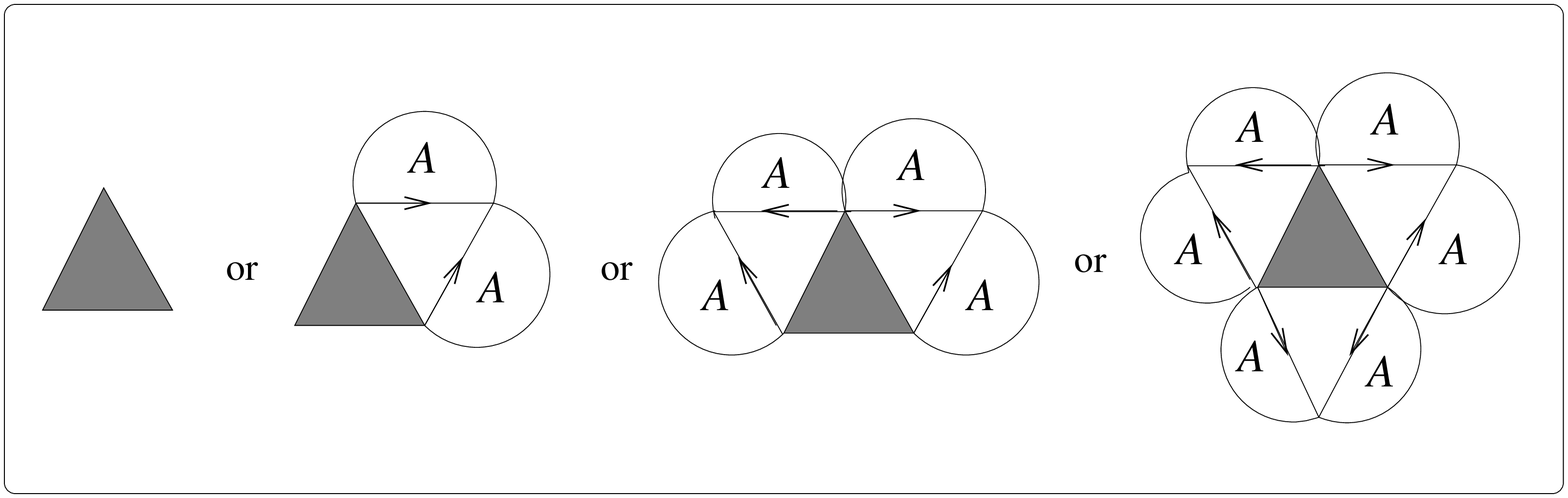}}
  \caption{ The species $\apt$}
  \label{new2}
\end{figure}

\begin{figure}[h]
 \centerline{\includegraphics[height=.18\textheight]{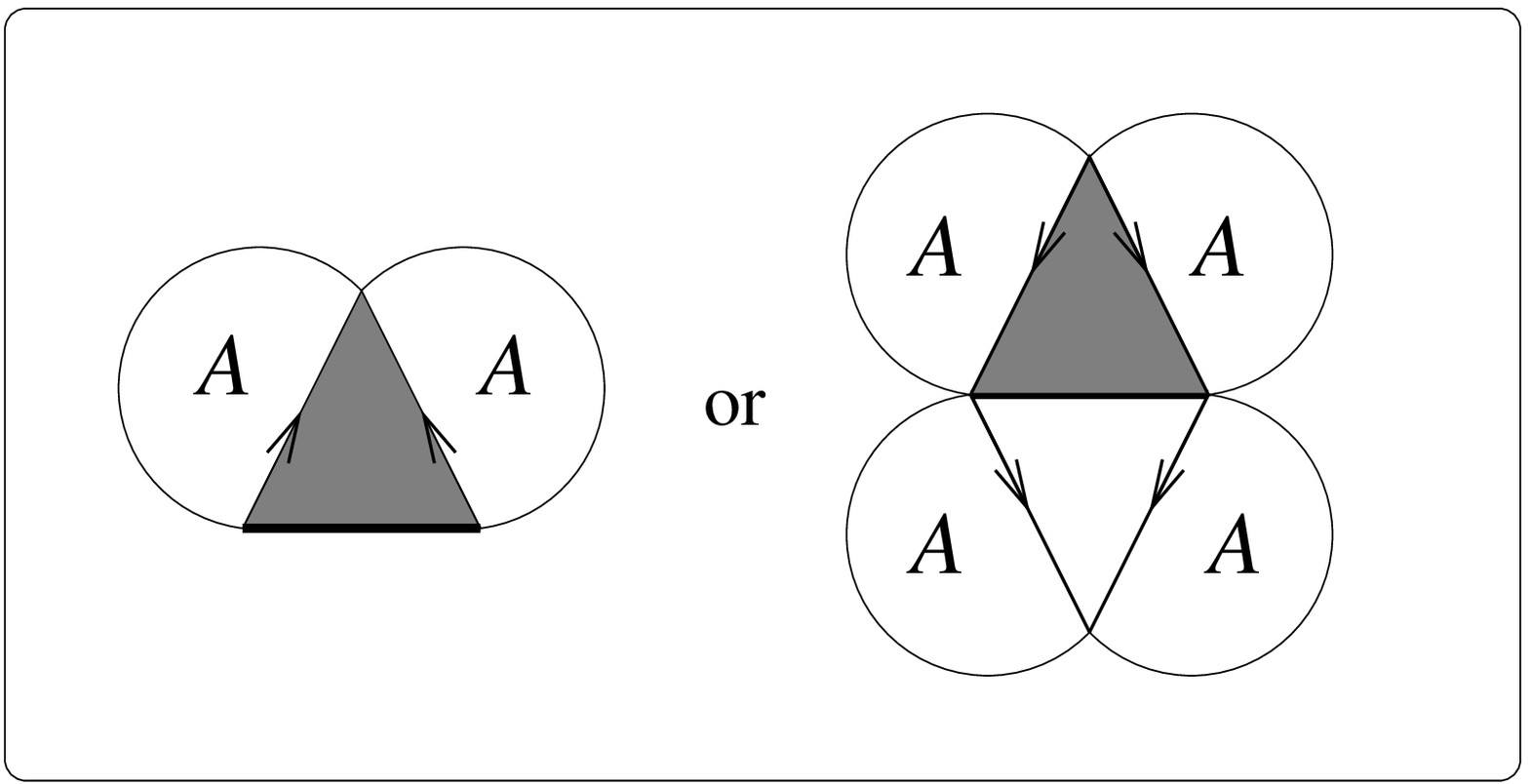}}
  \caption{The species $\apbt$}
  \label{new3}
\end{figure}

This subsection is devoted to planar \ars{}, \ie \ars{} admitting an embedding in the plane in such a way that all internal faces are triangles. The difference here is that the embedding is not explicitely given and that reflexive symmetries are possible. In other words, planar 2-trees are viewed as simple graphs. The dissymetry theorem for the species $\ap$ of planar 2-trees yields
\begin{equation}
\ap=\apb+\apt-\apbt.\label{DTHP}
\end{equation}
Moreover, we have the following expressions for the pointed species $\apb$, $\apt$ and $\apbt$, in terms of the auxiliary species $P_4^{\rm bic}(X,Y)$ and $P_6^{\rm bic}(X,Y)$ introduced in Section 2. 
\begin{theo}
The species of pointed planar 2-trees $\apb$, $\apt$ and $\apbt$ satisfy the following isomorphisms of species :
\begin{eqnarray}
\apb(X)&=&1+X E_2(A)+P_4^{\rm bic}(X,Y)\arrowvert_{Y:=A},\label{pla1}\\
\apt(X)&=&X+X^2 E_2(A)+XE_2(A_+)+XP_6^{\rm bic}(X,Y)\arrowvert_{Y:=A},\label{pla2}\\
\apbt(X)&=&XE_2(A)+X^2E_2(A^2).\label{pla3}
\end{eqnarray}
\end{theo}
\begin{preuve}
 We obtain the functional equations (\ref{pla1}) and (\ref{pla3}) by analyzing the structures according to the degree of the distinguished edge. For example, the three terms on the right hand side of (\ref{pla1}) correspond respectively to the degrees 0, 1 and 2 of the pointed edge. This isomorphism is described in Figure~\ref{new1}. In (\ref{pla2}), the four terms correspond to the four possibilities for the number of edges of degree 2 in the pointed triangle, from 0 to 3; see Figure~\ref{new2}. For (\ref{pla3}), see Figure~\ref{new3}.
\end{preuve}
%
%

Combining the molecular expansion of the quotient species $P_4^{\rm bic}(X,A)$ and $P_6^{\rm bic}(X,A)$ established in Section 3 with Proposition~\ref{prop1} and Proposition~\ref{lem3}, gives the molecular expansion of the species $\apb$ and $\apt$. Note that we use the same notation for the coefficients of the different molecular expansions in the four following theorems.
\begin{theo}
The molecular expansion of the species $\apb$ of edge pointed planar \ars{} is given by
\begin{eqnarray}
\apb(X)&=&1+\sum_{k\geq 0}a_{k}^{1}X^k+\sum_{k\geq 1}a_{k}^{2}E_2(X^k)+\sum_{k\geq 1}a_{k}^{3}XE_2(X^k) \nonumber\\
&&\mbox{}+
\sum_{n\geq 1}a_{k}^{4}X^2E_2(X^k)+\sum_{k\geq 1}a_{k}^{5}P_4^{\rm bic}(X,X^k),\label{decc1}
\end{eqnarray}
where 
\begin{eqnarray}
a_{k}^1&=&{1\over 4}\cat_{k+1}-{3\over 4}\cat_{k\over 2}-{1\over 2}\cat_{k-1\over 2}+{1\over 2}\cat_{k-2\over 4},\nonumber\\
a_{k}^2&=&\cat_{k}-\cat_{\frac{k-1}{2}},\nonumber\\
a_{k}^3&=&a_{k}^5\ =\  \cat_k,\\
a_{k}^4&=&\frac{1}{2}(\cat_{k+1}-\cat_{\frac{k}{2}}).\nonumber
\end{eqnarray}
\end{theo}
\begin{theo}
The molecular expansion of the species $\apt$ is given by
\begin{eqnarray}
\apt(X)&=&1+  \sum_{k\geq 0}a_{k}^{1}X^k+\sum_{k\geq 1}a_{k}^{2}X\cdot E_2(X^k)+\sum_{k\geq 2}a_{k}^{3}X^2E_2(X^k)\nonumber\\
&&\mbox{}+\sum_{k\geq 2}a_{k}^{4}XC_3(X^k)+\sum_{k\geq 2}a_{k}^{5}XP_6^{\rm bic}(X,X^k),\label{decc3}
\end{eqnarray}
where 
\begin{eqnarray}
 a_{k}^1&=&{1\over 6}(\cat_{k+1}-\cat_{k})-{1\over 2}\cat_{k\over 2}-\cat_{k-2\over 2}-{1\over 2}\cat_{k-1\over 2}-{1\over 6}\cat_{k-1\over 3}+{1\over 2}\cat_{k-4\over 6},\nonumber\\
a_{k}^2&=& a_{k}^5\ =\ \cat_k,\nonumber\\
a_{k}^3&=&\cat_{k+1}-\cat_{{k-1\over 3}},\\
a_{k}^4&=&\frac{1}{2}(\cat_k-\cat_{k-1\over 2}).\nonumber
\end{eqnarray}
\end{theo}
Proposition~\ref{prop1} and Proposition~\ref{lem3} also allow us to obtain the molecular expansion of the species $\apbt$.
\begin{theo}
The molecular expansion of the species $\apbt$ of planar \ars{} pointed at a triangle with a distinguished edge is given by
\begin{equation}
\apbt(X)=\sum_{k\geq 0}a_{k}^{1}X^k+\sum_{k\geq 1}a_{k}^{2}XE_2(X^k)+
\sum_{k\geq 1}a_{k}^{3}X^2E_2(X^k),\label{decc2}
\end{equation}
where 
\begin{eqnarray}
a_{k}^1&=&\frac{1}{2}\big(\cat_{k+1}-\cat_k-\cat_{{k-1\over 2}}-\cat_{{k\over 2}}\big),\nonumber\\
a_{k}^2&=&\cat_k,\\
a_{k}^3&=&\cat_{k+1}.\nonumber
\end{eqnarray}
\end{theo}

Using the dissymetry theorem, we are now able to put together relations (\ref{decc1})-(\ref{decc3})-(\ref{decc2}) and give an explicit form of the molecular expansion of the species $\ap$ of planar \ars{}.
\begin{theo}\label{Decmolplanaire}
The molecular expansion of the species $\ap$ of planar \ars{} 
is given by the following formula 
\begin{eqnarray}
\ap(X)&=&1+  \sum_{k\geq 1}a_{k}^{1}X^k+\sum_{k\geq 1}a_{k}^{2}E_2(X^k)+\sum_{k\geq 1}a_{k}^{3}XE_2(X^k)+
\sum_{k\geq 2}a_{k}^{4}X^2E_2(X^k)\nonumber\\
&&\mbox{}+\sum_{k\geq 2}a_{k}^{5}XC_3(X^k)+\sum_{k\geq 0}a_{k}^{6}P_4^{\rm bic}(X,X^k)+\sum_{k\geq 0}a_{k}^{7}XP_6^{\rm bic}(X,X^k),\label{dectot}
\end{eqnarray}
where 
\begin{eqnarray}
a_{k}^1&=&-{1\over 12}\cat_{k+1}+{1\over 3}\cat_{k}-{3\over 4}\cat_{k\over 2}-{1\over 2}\cat_{k-1\over 2}-{1\over 6}\cat_{k-1\over3}+{1\over 2}\cat_{k-2\over 4}+{1\over 2}\cat_{k-4\over 6},\nonumber \\
a_{k}^2&=&\cat_{k}-\cat_{\frac{k-1}{2}},\nonumber\\
a_{k}^3&=& a_{k}^6\ =\  a_{k}^7\ =\   \cat_k, \label{cc}     \\
a_{k}^4&=&{1\over 2}(\cat_{k+1}-\cat_{k\over 2})-\cat_{{k-1\over 3}},\nonumber\\a_{k}^5&=&{1\over 2}(\cat_k-\cat_{k-1\over 2}).\nonumber
\end{eqnarray}
\end{theo}


\section{Enumeration formulas}
\subsection{Enumeration of plane 2-trees}
Before obtaining the explicit enumeration of plane 2-trees, we recall some basic formulas involving index series of the species of 2-element sets ($E_2$) and of oriented 3-cycles ($C_3$) :
\begin{eqnarray}
Z_{E_2}(x_1,x_2,\ldots)&=&{1\over 2}(x_1^2+x_2),\quad \ \,  \Gamma_{E_2}(x_1,x_2,\ldots)={1\over 2}(x_1^2-x_2),\\
Z_{C_3}(x_1,x_2,\ldots)&=&{1\over 3}(x_1^3+2x_3),\quad \Gamma_{C_3}(x_1,x_2,\ldots)={1\over 3}(x_1^3-x_3).
\end{eqnarray}  

We will also use some substitutional laws of the theory of species : for any species $F$ and $G$ such that $G(0)=0$ ($G$ has no structure on the empty set), we have
\begin{eqnarray}
(F\circ G)(x)&=&F(G(x)),\label{sub1}\\
({F\circ G})^{\sim}(x)&=&Z_F(\widetilde{G}(x),\widetilde{G}(x^2),\ldots ).\label{sub2},\\
(\overline{F\circ G})(x)&=&\Gamma_F(\bar{G}(x),\bar{G}(x^2),\ldots ),\label{sub3}\\
Z_{F\circ G}&=&Z_{F}\circ Z_G,\label{sub4}\\
\Gamma_{F\circ G}&=&\Gamma_{F}\circ \Gamma_G,\label{sub5}
\end{eqnarray}
where $\circ$ denotes the plethystic composition on the right hand side of (\ref{sub4}) and (\ref{sub5}). 

If the species $G$ has some structures on the empty set, \ie $G(0)=g_0\neq 0$, formulas (\ref{sub2})--(\ref{sub4}) remain valid. However, formula (\ref{sub1}) should then be replaced by
\begin{equation}
(F\circ G)(x) = Z_{F}(G(x),g_0,g_0,\ldots ),
\end{equation}
and there is no known general formula for $\Gamma$. Here, we only need the following formulas
\begin{eqnarray}
\Gamma_{E_2(G)}(x_1,x_2,\ldots)&=&g_0+{1\over 2}(\Gamma_G^2(x_1,x_2,\ldots)-\Gamma_G(x_2,x_4,\ldots)),\\
\Gamma_{C_3(G)}(x_1,x_2,\ldots)&=&g_0+{1\over 3}(\Gamma_G^3(x_1,x_2,\ldots)-\Gamma_G(x_3,x_6,\ldots)).
\end{eqnarray}

We now give the explicit enumerative formulas provided directly by the molecular expansion of the species of plane 2-trees.
\begin{theo}
The numbers $a_{\pi,n}$, $\widetilde{a}_{\pi,n}$  and  $\overline{a}_{\pi,n}$ of labelled, unlabelled and unlabelled asymmetric plane \ars{} on $n$ triangles, $n\geq 2$, are given by
\begin{eqnarray}
a_{\pi ,n}&=&n!({2\over 3}\cat_n-{1\over 6}\cat_{n+1}),\\
{\widetilde a}_{\pi ,n}&=&{2\over 3}\cat_n-{1\over 6}\cat_{n+1}+{1\over 2}\cat_{{n\over 2}}+{2\over 3}\cat_{{n-1\over 3}}\label{asym},\\
\overline{a}_{\pi,n} &=&{2\over 3}\cat_n-{1\over 6}\cat_{n+1}-{1\over 2}\cat_{k\over 3}-{1\over 3}\cat_{n-1\over 3}.
\end{eqnarray}
\end{theo}

To obtain these enumerating formulas, we can also use the expressions (\ref{t})--(\ref{ttt}) which lead to closed formulas for the associated series of the three pointed species : the exponential generating series of labelled structures,
\begin{eqnarray}
\abb_{\pi}(x)&=&{1\over 2}(1+A^2(x)),\nonumber\\
\abt_{\pi}(x)&=&{x\over 3}(2+A^3(x)),\label{l1}\\
\abbt_{\pi}(x)&=&A^2(x)-A(x),\nonumber
\end{eqnarray}
the ordinary generating series of unlabelled structures
\begin{eqnarray}
\widetilde{\ab}_{\pi }^{-}(x)&=&{1\over 2}(A^2(x)+A(x^2)),\nonumber\\
\widetilde{\ab}_{\pi }^{\vartriangle}(x)&=&{x\over 3}(A^3(x)+2A(x^3)),\label{u1}\\
\widetilde{\ab}_{\pi }^{\:\underline{\!\vartriangle\!}}(x)&=&A^2(x)-A(x),\nonumber
\end{eqnarray}
the cycle index series
\begin{eqnarray}
Z_{\abb_{\pi}}(x_1,x_2,\ldots )&=&{1\over 2}\big(A^2(x_1)+A(x_2)\big),\nonumber\\
Z_{\abt_{\pi}}(x_1,x_2,\ldots )&=&{x_1\over 3}\big(A^3(x_1)+2A(x_3)\big),\label{z1}\\
Z_{\abbt_{\pi}}(x_1,x_2,\ldots )&=&A^2(x_1)-A(x_1),\nonumber
\end{eqnarray}
the asymmetry cycle index series
\begin{eqnarray}
\Gamma_{\abb_{\pi}}(x_1,x_2,\ldots )&=&1+{1\over 2}\big(A^2(x_1)-A(x_2)\big),\nonumber\\
\Gamma_{\abt_{\pi}}(x_1,x_2,\ldots )&=&x_1+{x_1\over 3}\big(A^3(x_1)-A(x_3)\big),\label{g1}\\
\Gamma_{\abbt_{\pi}}(x_1,x_2,\ldots )&=&A^2(x_1)-A(x_1).\nonumber
\end{eqnarray}

We emphasize the fact, used above, that since the species $A$ is asymmetric we have the following relations
\begin{equation}
A(x)=\widetilde{A}(x)\ =\ \overline{A}(x)\quad \mbox{and}\quad Z_A(x_1,x_2,\ldots )=A(x_1)=\Gamma_A(x_1,x_2,\ldots).
\end{equation}
We then deduce easily (thanks to the dissymmetry theorem) the expressions of the series associated with the species of plane 2-trees
\begin{prop}
The series associated to the species $\ab_{\pi}$ of plane 2-trees are given by 
\begin{eqnarray}
\ab_{\pi}(x)&=&{1\over 2}+{2\over 3}x+A(x)-{1\over 2}A^2(x)+{x\over 3}A^3(x),\nonumber   \\
\widetilde{\ab}_{\pi}(x)&=&1+x+A(x)+{x\over 3}A^3(x)-{1\over 2}A(x^2)-{x\over 3}A(x^3)-{1\over 2}A^2(x),\nonumber   \\
\bar{\ab}_{\pi}(x)&=&A(x)+{x\over 3}A^3(x)-A(x^2)-A(x^3)-{1\over 2}A^2(x),\\
Z_{\ab_{\pi}}(x_1,x_2,\ldots )&=&A(x_1)+{1\over 2}A(x_2)+{2\over 3}x_1A(x_3)-{1\over 2}A^2(x_1)+{x_1\over 3}A^3(x_1),\nonumber   \\
\Gamma_{\ab_{\pi}}(x_1,x_2,\ldots)&=&1+x_1+A(x_1)+{x_1\over 3}A^3(x_1)-{1\over 2}A(x_2)-{x_1\over 3}A(x_3)-{1\over 2}A^2(x_1). \nonumber   
\end{eqnarray}
\end{prop}
To recover the formulas (\ref{asym}), we can use the dissymmetry theorem and the next proposition giving the enumeration of the different pointed plane 2-trees.
\begin{prop}
The coefficients $a_{\pi,n}^{-}$, $a_{\pi,n}^{\vartriangle}$, $a_{\pi,n}^{\: \underline{\!\vartriangle\!}}$ representing the numbers of labelled structures with $n$ triangles for the different pointings, $\widetilde{a}_{\pi,n}^{-}$, $\widetilde{a}_{\pi,n}^{\vartriangle}$, $\widetilde{a}_{\pi,n}^{\underline{\!\vartriangle\!}}$ for the numbers of unlabelled structures, and  $\overline{a}_{\pi,n}^{-}$, $\overline{a}_{\pi,n}^{\vartriangle}$, $\overline{a}_{\pi,n}^{\underline{\!\vartriangle\!}}$ for unlabelled asymmetric structures, are given, for $n\geq 2$, by
\begin{eqnarray}
a_{\pi ,n}^{-}&=&{n!\over 2}\cat_{n+1},\qquad \quad \qquad \quad \ \ \nonumber   \\
a_{\pi ,n}^{\vartriangle}&=&{n!\over 3}(\cat_{n+1}-\cat_n),\label{zb2}\\
a_{\pi ,n}^{\underline{\!\vartriangle\!}}&=&n!(\cat_{n+1}-\cat_n),\nonumber   
\end{eqnarray}
\begin{eqnarray}
\widetilde{a}_{\pi ,n}^{-}&=&{1\over 2}(\cat_{n+1}+\cat_{{n\over 2}}),\nonumber   \\
\widetilde{a}_{\pi ,n}^{\vartriangle}&=&{1\over 3}(\cat_{n+1}-\cat_n+2\cat_{{n-1\over3}}),\\
\widetilde{a}_{\pi ,n}^{\underline{\!\vartriangle\!}}&=&\cat_{n+1}-\cat_n,\nonumber   
\end{eqnarray}
and
\begin{eqnarray}
\bar{a}_{\pi ,n}^{-}&=&{1\over 2}(\cat_{n+1}-\cat_{{n\over 2}}),\nonumber   \\
\bar{a}_{\pi ,n}^{\vartriangle}&=&{1\over 3}(\cat_{n+1}-\cat_n-\cat_{{n-1\over3}}),\label{asym-}\\
\bar{a}_{\pi ,n}^{\underline{\!\vartriangle\!}}&=&\cat_{n+1}-\cat_n.\nonumber   
\end{eqnarray}
\end{prop}

\begin{preuve}
To obtain these coefficients, we  simply use relations (\ref{l1}), (\ref{u1}) and (\ref{g1}).
\end{preuve}

We now give the explicit expressions for the cycle index series of the species of plane \ars{}.
\begin{prop}
The cycle index series and the asymmetric index series of the species of plane \ars{} are
\begin{equation}
Z_{\ab_{\pi}}(x_1,x_2,\ldots )=1+\sum_{n\geq 1}({2\over 3}\cat_{n}-{1\over 6}\cat_{n+1})x_1^n+{1\over 2}\sum_{n\geq 1}\cat_nx_2^n+{2\over 3}x_1\sum_{n\geq 0}\cat_nx_3^n,
\end{equation}
\begin{equation}
\Gamma_{\ab_{\pi}}(x_1,x_2,\ldots )=1+x_1+\sum_{n\geq 1}({2\over 3}\cat_{n}-{1\over 6}\cat_{n+1})x_1^n-{1\over 2}\sum_{n\geq 1}\cat_nx_2^n-{1\over 3}x_1\sum_{n\geq 0}\cat_nx_3^n.
\end{equation}
\end{prop}
\begin{preuve}
We first express the cycle index series given by the relations (\ref{z1}) in powers of $x_1$, $x_2$, \ldots
\begin{eqnarray}
Z_{\abb_{\pi}}(x_1,x_2,\ldots )&=&{1\over 2}\sum_{n\geq 0}\cat_{n+1}x_1^n+{1\over 2}\sum_{n\geq 0}\cat_nx_2^n,\nonumber\\
Z_{\abt_{\pi}}(x_1,x_2,\ldots )&=&{1\over 3}\sum_{n\geq 1}(\cat_{n+1}-\cat_n)x_1^n+{2\over 3}x_1\sum_{n\geq 0}\cat_nx_3^n,\\
Z_{\abbt_{\pi}}(x_1,x_2,\ldots )&=&\sum_{n\geq 1}(\cat_{n+1}-\cat_n)x_1^n.\nonumber
\end{eqnarray}
We also have
\begin{eqnarray}
\Gamma_{\abb_{\pi}}(x_1,x_2,\ldots )&=&1+{1\over 2}\sum_{n\geq 0}\cat_{n+1}x_1^n-{1\over 2}\sum_{n\geq 0}\cat_nx_2^n\nonumber,\\
\Gamma_{\abt_{\pi}}(x_1,x_2,\ldots )&=&x_1+{1\over 3}\sum_{n\geq 1}(\cat_{n+1}-\cat_n)x_1^n-{1\over 3}x_1\sum_{n\geq 0}\cat_nx_3^n,\\
\Gamma_{\abbt_{\pi}}(x_1,x_2,\ldots )&=&\sum_{n\geq 1}(\cat_{n+1}-\cat_n)x_1^n.\nonumber
\end{eqnarray}
It suffices then to use the dissymmetry theorem to obtain the stated result.
\end{preuve}

\subsection{Enumeration of planar 2-trees}

We now give all associated series of the species $\apb$, $\apt$ and $\apbt$ using substitutionnal laws of the theory of species. After this, we will be able to give all coefficients arising in these differents series, and, with the dissymetry theorem, we obtain the number of labelled and unlabelled planar \ars{} on $n$ triangles as well as the coefficients of its cycle and asymmetry index series.
\begin{theo}\label{th14}
The exponential generating function of labelled stuctures for the species $\apb$, $\apt$ and $\apbt$ of planar pointed \ars{}
are given, in terms of the species $A$, by
\begin{eqnarray}
\apb(x)&=&1+{x\over 2}(1+A^2(x))+{1\over 4}x^2A^4(x),\nonumber\\
\apt(x)&=&x+{x^2\over 2}(1+A^2(x))+{x\over 2}A_+^2(x)+{x^4\over 6}A^6(x),label{gf1}\\
\apbt(x)&=&{x\over 2}(1+A^2(x))+{x^2\over 2}(1+A^4(x)).\nonumber
\end{eqnarray}
Moreover, the ordinary generating series of unlabelled structures of these species are given by
\begin{eqnarray}
\widetilde{\ab}_{\mathrm{p}}^{-}(x)
&\makebox[0mm][c]{=}&
1+xA(x)+{x\over 2}(A^2(x)+A(x^2))+{x^2\over 4}(A^4(x)+3A^2(x^2))\,,\nonumber\\
\widetilde{\ab}_{\mathrm{p}}^{\tr}(x)
&\makebox[0mm][c]{=}&
x+{x^2\over 2}(A^2(x)+A(x^2))+{x\over 2}(A_+^2(x)+A_+(x^2))\nonumber\\
&&\hspace*{4cm}+{x^4\over 6}(A^6(x)+2A^2(x^3)+3A^3(x^2))\,,\\
\widetilde{\ab}_{\mathrm{p}}^{\btr}(x)
&\makebox[0mm][c]{=}&
{x\over 2}(A^2(x)+A(x^2))+{x^2\over 2}(A^4(x)+A^2(x^2)).\nonumber
\end{eqnarray}
\end{theo}
\begin{cor}
The exponential and the ordinary generating functions of the species of planar 2-trees are given, in terms of $A$, by 
\begin{eqnarray}
\ap(x)&=&1+x+{x\over 2}A_+^2(x)+{x^2\over 2}A^2(x)-{x^2\over 4}A^4(x)-{x^4\over 6}A^6(x),\nonumber\\ 
\widetilde{\ab}_{\mathrm{p}}(x)&=&1+x+{x\over 2}(A_+^2(x)+A_+(x^2))+{x^2\over 2}A(x^2)+{x^2\over 2}(A^2(x)-A^2(x^2))\\
&&\hspace*{4cm}-{x^2\over 4}A^4(x)+{x^4\over 6}(A^6(x)+2A^2(x^3)+3A^3(x^2)).\nonumber
\end{eqnarray}
\end{cor}

A simple extraction of coefficients in Theorem~\ref{th14}, combined with Proposition~\ref{prop1}, yields the following corollary.
\begin{cor}
The numbers $a_{\mathrm{p},n}^{-}$, $a_{\mathrm{p},n}^{\vartriangle} $ and $a_{\mathrm{p},n}^{\: \underline{\!\vartriangle\!}} $ of labelled planar \ars{} on $n$ triangles pointed respectively at an edge, at a triangle, and at a triangle pointed at one of its edges, are given by
\begin{eqnarray}
a_{\mathrm{p},n}^{-}&=&{n!\over 4}\cat_{n+1},\nonumber\\
a_{\mathrm{p},n}^{\vartriangle}&=&{n!\over 6}(\cat_{n+1}-\cat_n) ,\\
a_{\mathrm{p},n}^{\: \underline{\!\vartriangle\!}}&=&{n!\over 2}(\cat_{n+1}-\cat_n).\nonumber
\end{eqnarray}
Moreover, for the same pointed series, the numbers of unlabelled structures on $n$ triangles $\widetilde{a}_{\mathrm{p},n}^{-}$, $\widetilde{a}_{\mathrm{p},n}^{\vartriangle}$ and $\widetilde{a}_{\mathrm{p},n}^{\: \underline{\!\vartriangle\!}}$ have the following expressions :
\begin{eqnarray}
\widetilde{a}_{\mathrm{p},n}^{-}&=&{1\over 4}\cat_{n+1}+{1\over 2}\cat_{{n-1\over 2}}+{3\over 4}\cat_{{n\over 2}},\nonumber\\
\widetilde{a}_{\mathrm{p},1}^{\vartriangle}&=&1,\ \widetilde{a}_{\mathrm{p},2}^{\vartriangle}=1,\nonumber\  \widetilde{a}_{\mathrm{p},3}^{\vartriangle}=2, \  \widetilde{a}_{\mathrm{p},4}^{\vartriangle}=6,\\
\widetilde{a}_{\mathrm{p},n}^{\vartriangle}&=&{1\over 6}(\cat_{n+1}-\cat_n)+{1\over 2}\big(\cat_{{n-1\over2}}+\cat_{{n\over 2}})\big)+{1\over 3}\cat_{{n-1\over 3}},\quad n\geq 5\\
\widetilde{a}_{\mathrm{p},n}^{\: \underline{\!\vartriangle\!}}
&=&{1\over 2}(\cat_{n+1}-\cat_n)+\cat_{{n-1\over 2}}+\cat_{{n\over2}}.\nonumber
\end{eqnarray}
\end{cor}
Hence, the dissymmetry theorem leads us to enumeration formulas for labelled and unlabelled planar \ars{} as follows. For the unlabelled asymmetric enumeration, we use directly the molecular decomposition of the species $\ap$.
\begin{theo}\label{unlabelled}
The numbers $a_{\mathrm{p},n}$, $\widetilde{a}_{\mathrm{p},n}$ and $\bar{a}_{\mathrm{p},n}$ of labelled, unlabelled and unlabelled asymmetric planar \ars{} on $n$ triangles, are given by the following formulas
\begin{eqnarray}
a_{\mathrm{p},n}&=&n!({1\over 3}\cat_n-{1\over 12}\cat_{n+1}),\\
\widetilde{a}_{\mathrm{p},n}&=&{1\over3}\cat_n-{1\over 12}\cat_{n+1}+{1\over 2}\cat_{{n-1\over 2}}+{1\over 3}\cat_{{n-1\over 3}}+{3\over 4}\cat_{{n\over2}},\label{der}\\
\bar{a}_{\mathrm{p},n}&=& -{1\over 12}\cat_{n+1}+{1\over 3}\cat_{n}-{3\over 4}\cat_{n\over 2}-{1\over 2}\cat_{n-1\over 2}-{1\over 6}\cat_{n-1\over3}+{1\over 2}\cat_{n-2\over 4}+{1\over 2}\cat_{n-4\over 6}.
\end{eqnarray}
\end{theo}
Finally, we give the expression of the asymmetry index series of the species $\ap$ of planar 2-trees obtained directly from the molecular expansion of the species $\ap$.
\begin{prop}
The asymmetry index series of the species of planar 2-trees is given by
$$
\Gamma_{\ap}(x_1,x_2,\ldots)=1+x_1+\sum_{n}\gamma_{n}^{1}x_1^n+\sum_{n}\gamma_n^{2}x_2^n+\sum_{n}\gamma_{n}^{3}x_1x_2^n+\sum_{n}\gamma_{n}^{4}x_1^2x_2^n+
$$
\begin{equation}
+\sum_{n}\gamma_n^{5}x_1x_3^n+\sum_{n}\gamma_{n}^{6} x_{2}x_{4}^{n}+\sum_{n}\gamma_n^{7}x_1x_3x_6^n,
\end{equation}
where
\begin{eqnarray}
\gamma_n^{1} & = &-{1\over 12}\cat_{n+1}+{1\over 3}\cat_n ,\nonumber\\
\gamma_n^{2} & = & \gamma_n^{3} \ =\ -{1\over 2}\cat_n,\nonumber\\
\gamma_n^{4} & = & -{1\over 4}\cat_{n+1},\\
\gamma_n^{5} & = & -{1\over 6}\cat_n,\nonumber\\
\gamma_n^{6} & = & \gamma_{n}^{7} \ = \ {1\over 2}\cat_n.\nonumber
\end{eqnarray}
\end{prop}

\subsection{Another method for the unlabelled enumeration}
In order to obtain the unlabelled enumeration of plane and planar 2-trees, we can also use the approach of Palmer and Read in \cite{PR}. Remark first that, for any species $F$, we can write
\begin{equation}
F=\sum_{k\geq 1}F_{(k)},
\end{equation}
where for $k\geq 1$, $F_{(k)}$ represents the symmetric part of $F$  of order $k$, \ie the subspecies consisting of $F$-structures whose stabilizer is of order $k$ exactly. In particular, $F_{(1)}=\overline{F}$, the asymmetric part of $F$.\\

Also note that, for $G=F_{(k)}$, $k\geq 1$, we have $G(x)={1\over k}\widetilde{G}(x)$, since an unlabelled $F_{(k)}$-structure of degree $n$ can be labelled in $n!/k$ ways. Hence
\begin{equation}
\widetilde{F}(x) = F(x)+\sum_{k\geq 2}{k-1\over k}\widetilde{F}_{(k)}(x).
\end{equation}
For plane 2-trees, we have

\begin{equation}
\ab_{\pi}=\overline{\ab}_{\pi}+\ab_{\pi,(2)}+\ab_{\pi,(3)},
\end{equation}
and for planar 2-trees,
\begin{equation}
\ap=\overline{\ab}_{\mathrm{p}}+  \ab_{{\rm p},(2)}+\ab_{{\rm p},(3)}+\ab_{{\rm p},(4)}+\ab_{{\rm p},(6)}.
\end{equation}
Hence, we can write
\begin{equation}
\widetilde{\ab}_{\pi}(x)=\ab_{\pi}(x)+{1\over2}\widetilde{\ab}_{\pi,(2)}(x)+{2\over3}\widetilde{\ab}_{\pi,(3)}(x),\label{PREAD1}
\end{equation}
and 
\begin{equation}\label{PREAD2}
\widetilde{\ab}_{\mathrm{p}}(x)=\ab_{\mathrm{p}}(x)+{1\over 2}\widetilde{\ab}_{\mathrm{p},(2)}(x)+{2\over 3}\widetilde{\ab}_{\mathrm{p},(3)}(x)+{3\over 4}\widetilde{\ab}_{\mathrm{p},(4)}(x)+{5\over 6}\widetilde{\ab}_{\mathrm{p},(6)}(x).
\end{equation} 
After identifying all terms appearing in (\ref{PREAD1}), we then deduce
\begin{equation}
\widetilde{\ab}_{\pi}(x)=\ab_{\pi}(x)+{1\over2}A(x^2)+{2\over3}xA(x^3),\label{palm1}\\
\end{equation}
for the plane case. For planar 2-trees, we have
\begin{eqnarray}
\widetilde{\ab}_{{\rm p},(2)} (x) &=& {3\over 2}(A(x^2)-x^2A(x^4))+xA(x^2)-x^4A(x^6),\nonumber\\
\widetilde{\ab}_{{\rm p},(3)} (x) &=& {1\over 2}(xA(x^3)-x^4A(x^6)),\\
\widetilde{\ab}_{{\rm p},(4)} (x) &=& x^2A(x^4),\quad \widetilde{\ab}_{{\rm p},(6)} \ =\ x^4A(x^6),\nonumber
\end{eqnarray}
which yields
\begin{equation}\label{palm2}
 \widetilde{\ab}_{\mathrm{p}}(x)=\ab_{\mathrm{p}}(x)+{1\over 2}xA(x^2)+{1\over 3}xA(x^3)+{3\over 4}A(x^2).
\end{equation}
It remains to extract the  coefficients of $x^n$ in equations (\ref{palm1}) and (\ref{palm2}) to find the numbers of unlabelled plane and planar 2-trees over $n$ triangles, given by (\ref{asym}) and (\ref{der}).

{\it E-mail address} : {\bf \{gilbert, lamathe, leroux\}@lacim.uqam.ca}\\
\\
* Corresponding author : Pierre Leroux\\ 
{LaCIM, D\'epartement de Math\'ematiques, UQ\`AM\\
C. P. 8888, succursale Centre-Ville\\
Montr\'eal (Qc) Canada H3C 3P8 


\begin{thebibliography}{33}
\bibitem{ALL1} P. Auger, G. Labelle, and P. Leroux, {\it Combinatorial addition formulas}, Proceedings FPSAC'01, Tempe, Arizona, May 21-25 2001, H. Barcelo and V. Welker, Eds, pp 19--26.

\bibitem{ALL2} P. Auger, G. Labelle, and P. Leroux, {\it Combinatorial addition formulas and applications}, Advances in Applied Mathematics, to appear.

\bibitem{BLL} F. Bergeron, G. Labelle, and P. Leroux, {\it Combinatorial Species and tree-like structures}, Encyclopedia of Mathematics and it'ss Applications, vol. 67, Cambridge University Press, (1998).


\bibitem{BR} W. G. Brown, {\it Enumeration of triangulations of the disk}, Proc. London Math. Soc. {\bf 14}, 746-768, (1964).

\bibitem{CBBCL} S. J. Cyvin, J. Brunvoll, E. Brensdal, B. N. Cyvin and E. K. Lloyd, {\it Enumeration of Polyene Hydrocarbons : A Complete Mathematical Solution}, J. Chem. Inf. Comput. Sci., {\bf 35}, 743-751, (1995).

\bibitem{TF1} T. Fowler, I. Gessel, G. Labelle, P. Leroux, {\it Specifying \ars{}}, Proceedings FPSAC'00, Moscou, 26-30 juin 2000, 202-213.

\bibitem{TF2} T. Fowler, I. Gessel, G. Labelle, P. Leroux, {\it The Specification of \ars{}}, Advances in Applied Mathematics, to appear.

\bibitem{HP} F. Harary and E. Palmer, {\it Graphical Enumeration}, Academic Press, New York, (1973).

\bibitem{HPR} F. Harary, E. Palmer and R. Read, {\it On the cell-growth problem for arbitrary polygons}, Discrete Mathematics, 11, 371--389, (1975).

\bibitem{K} A. Kerber, {\it Enumeration under Finite Group Action: Symmetry Classes of Mappings},  Combinatoire \'enum\'erative, Proceedings, Montr\'eal, Qu\'ebec, Lectures Notes in Mathematics, vol. 1234, Springer-Verlag, New-York/Berlin, 160--176, (1985).

\bibitem{GL} G. Labelle, {\it On Asymmetric Structures}, Discrete Mathematics, 99, 141-162, (1992).

\bibitem{LLP} G. Labelle, J. Labelle and K. Pineau, {\it Sur une g\'en\'eralisation des s\'eries indicatrices d'esp\`eces}, J. of Comb. Theory, Series A, {\bf 69}, No. 1, 17--35, (1995).  

\bibitem{LLL} G. Labelle, C. Lamathe and P. Leroux, {\it D\'eveloppement mol\'eculaire de l'esp\`ece des 2-arbres planaires}, Proceedings GASCom 2001, 41--46, (2001).


\bibitem{JL} J. Labelle, {\it Quelques esp\`eces sur les ensembles de petite cardinalit\'e}, Annales des Sciences Math\'ematiques du Qu\'ebec, {\bf 11}, 31-58, (1985). 

\bibitem{PR} E. Palmer and R. Read, {\it On the Number of Plane \ars{}}, J. London Mathematical Society {\bf 6}, 583-592, (1973).

\bibitem{SP} N. J. A. Sloane and S. Plouffe, {\it The Encyclopedia of Integer Sequences}, Academic Press, San Diego, (1995).


\end{thebibliography}
\end{document}